\documentclass[11pt]{article}   
\usepackage{amssymb,amscd,latexsym}   
\usepackage{amsmath}
\usepackage{amsthm}
\usepackage[dvipsnames]{xcolor}
\usepackage{pstricks,pst-grad}
\newrgbcolor{l}{0.8 0.3 0.7}
\newrgbcolor{li}{0 0.3 0.7}
\newrgbcolor{lig}{1 0 1}
\newrgbcolor{ligh}{0 0.3 0}
\newrgbcolor{r}{1 0 0}
\usepackage{tikz}
\definecolor{grey}{rgb}{0.9,0.90,0.90}

\newcommand{\rar}{\rightarrow}
\newcommand{\lar}{\longrightarrow}
\newcommand{\llar}{-\kern-5pt-\kern-5pt\longrightarrow}
\newcommand{\surjects}{\twoheadrightarrow}

\newtheorem{Theorem}{Theorem}[section]
\newtheorem{Lemma}[Theorem]{Lemma}
\newtheorem{Corollary}[Theorem]{Corollary}
\newtheorem{Proposition}[Theorem]{Proposition}
\newtheorem{Remark}[Theorem]{Remark}
\newtheorem{Example}[Theorem]{Example}
\newtheorem{Conjecture}[Theorem]{Conjecture}

\newtheorem{Question}[Theorem]{Question}

\DeclareMathOperator{\depth}{depth}
\DeclareMathOperator{\Ass}{Ass}

\def\restr{{\kern-1pt\restriction\kern-1pt}}

\def\sqr#1#2{{\vcenter{\hrule height.#2pt
        \hbox{\vrule width.#2pt height#1pt \kern#1pt
            \vrule width.#2pt}
        \hrule height.#2pt}}}
\def\phi{\varphi}
\def\demo{\noindent{\bf Proof. }}
\def\square{\mathchoice\sqr64\sqr64\sqr{4}3\sqr{3}3}
\def\qed{\hspace*{\fill} $\square$}

\def\xx{{\bf x}}
\def\yy{{\bf y}}

\def\ff{{\bf f}}

\def\ff{{\bf f}}

\def\fm{{\mathfrak m}}

\def\hht{{\rm ht}\,}
\def\depth{{\rm depth}\,}

\def\rk{{\rm rank}\,}

\def\spec#1{{\rm Spec}\, (#1)}

\def\pp{{\mathbb P}}
\def\p{{\mathbb P}}

 2
 3
 4
 5
 4
 3
 2
 1
 1
 2
 3
 1
 0
 1
 2
 3
 4
 5

\let\oldsqrt\sqrt
\def\sqrt{\mathpalette\DHLhksqrt}
\def\DHLhksqrt#1#2{%
\setbox0=\hbox{$#1\oldsqrt{#2\,}$}\dimen0=\ht0
\advance\dimen0-0.2\ht0
\setbox2=\hbox{\vrule height\ht0 depth -\dimen0}%
{\box0\lower0.4pt\box2}}

\begin{document}

\begin{center}
{\Large{\bf\sc Homaloidal determinants}}\footnotetext{2010 AMS {\it Mathematics Subject
Classification}: 13B25, 13D02, 13H10, 13H15, 14E05, 14E07, 14J70, 14M05, 14M12.}

\vspace{0.3in}

{\large\sc Maral Mostafazadehfard}\footnote[1]{Under a
CNPq Doctoral Fellowship}
\quad\quad {\large\sc Aron  Simis}\footnote[2]{Partially
supported by a CNPq grant and a PVNS Fellowship from CAPES.}

\end{center}


\bigskip

\begin{abstract}

A form in a polynomial ring over a field is said to be homaloidal if its polar map is a Cremona map, i.e., if the rational map defined by the partial derivatives of the form has an inverse rational map.
The object of this work is the search for homaloidal polynomials that are the determinants of sufficiently structured matrices. We focus on generic catatalecticants, with special emphasis on the Hankel matrix. An additional focus is on certain degenerations or specializations thereof. In addition to studying the homaloidal nature of these determinants, one establishes several results on the ideal theoretic invariants of the respective gradient ideals, such as primary components, multiplicity, reductions and free resolutions.

\end{abstract}

\section{Introduction}

The subject of Cremona transformations  is a classical chapter of algebraic geometry. However, the
Cremona group of the projective space $\pp^ n$ is well understood only for $n\leq 2$. In higher dimension the structure of this group is far from being clarified, and the classification of such maps, except for a few special classes, is poorly known.
An important class of Cremona maps of  $\pp^ n$ comes off the so-called {\em polar maps}, that is, rational maps whose coordinates are the partial derivatives of a homogeneous polynomial $f$ in the homogeneous coordinate ring $R:=k[x_0,\ldots,x_n]$ of $\pp^n$.
Geometrically, the relevance of such a map is that
its indeterminacy locus is the singular locus of the corresponding hypersurface $V(f)$.

A homogeneous polynomial $f\in R$ whose polar map is a Cremona map is called {\em homaloidal} -- though more often this designation applies to the corresponding hypersurface rather than to $f$ itself.
Unfortunately, there are scarcely any general methods to studying, much less recognizing, such polynomials.
Results such as those in \cite{EKP} and \cite{KS}, although fascinating, are difficult to apply in practice and do not give a large picture.
A more circumscribed environment consists of polynomials that are determinants of square matrices with homogeneous entries of the same degree.
Alas, even for this class one lacks general methods, often happening that each such matrix requires a particular approach depending on its generic properties.

One of the goals of this paper is to consider  structured matrices whose entries are variables of a polynomial ring over a field.
Even for those there seems to be no comprehensive study of the homaloidal behavior of the corresponding determinants.
Still, one advantage of dealing with these matrices is that they are often $1$-generic in the sense
of \cite[Definition-Proposition 1.1]{Eisenbud2}. This implies that their determinant is irreducible
at the outset, thus allowing for a substantial class to search within.

Our approach is algebraic throughout.
Needless to justify, we will assume throughout that the base field has characteristic zero. A good deal of intuition about the results, if not the results themselves, gets lost in prime characteristic.
Anyway, from the geometric point of view, the study of a polar map in characteristic zero drives primevally through the properties of the Hessian determinant $h(f)$ of $f$, the reason being the classically known criterion for the dominance of the polar map in terms of the non
vanishing of the corresponding Hessian determinant. Although this criterion admits a vast generalization to arbitrary rational maps in terms of the Jacobian determinant of a basis of the corresponding linear system, it is in the polar case that the notion takes full role.

Though the overall objective is to detect homaloidal determinants and their properties, we soon became aware of a richness of notions from commutative algebra that come alongside in a natural way.
Often these notions and their use are crucial events for the geometric consequence. Thus, ideal theory in this paper is not just an aside, it is rather a live vein of the results.
For example, a major underlying problem in this context is to understand the properties of the so-called {\em gradient ideal} (or {\em singular ideal}) of the polynomial $f$,
that is, the ideal $J=J(f)\subset R$ generated by the partial derivatives of $f$.
Regardless of whether $f$ is homaloidal, a particular question asks
when $J$ has an irrelevant primary component, i.e., when it is not saturated. Of
course, if Proj$(R/(f)$ is smooth its gradient ideal will itself be an
irrelevant primary ideal (in addition, $f$ will not be homaloidal unless $\deg(f)\leq 2$).  One can easily cook up families of irreducible plane curves
whose gradient ideals have an irrelevant component. It is much harder to exhibit an
irreducible polynomial in $f\in k[x_0,x_1,x_2]$  such that its
gradient ideal has no irrelevant component (i.e., such that its gradient ideal is
perfect of codimension $2$). This question would naturally drive us in the theory of free divisors
(see \cite{LisbonProcDekker}, \cite{gradientsolo}, \cite{A-S}) which is slightly offshore our intention in this paper.
Nevertheless, the question itself is relevant in the present context as well.

The existence of tight lower bounds for the degrees of  syzygies of the gradient ideal has been around
in recent past, in connection with the so called inverse Poincar\'e problem  -- see, e.g., \cite{BruMen}, where it is shown that there are no syzygies in degree less than $\deg(f)-n$ provided the
singularities of $f$ are all normal crossings. Thus, in the context of this sort of singularities, only hypersurfaces $f$ of degree $\leq n+1$ qualify in order that the corresponding gradient ideal have any
linear syzygies at all.
This relates to the following problem which had been around as
part of the folklore of homaloidal polynomials.

\begin{Question}\label{degreebound}\rm Let $f\in k[\xx]=k[x_0,\ldots,x_n]$ denote a
squarefree form. If $f$ is homaloidal when is $\deg(f)\leq n+1$?
\end{Question}

It has been shown in \cite{CRS} that there are irreducible homaloidal forms
of degree arbitrarily larger than the number $n\geq 3$ of variables.
The classes of examples found are the dual hypersurfaces of certain rational scroll surfaces.
More recently, it has been announced in \cite[Section 1.3]{Huh2} the existence of many such examples rooted in Newton polytope theory as developed in \cite{Huh1}.
However, for certain  structured classes of homaloidal hypersurfaces
it is reasonable to expect that the degree stays beneath the number of variables.
Clearly, the notion of ``structured'' is quite foggy.
One sort of modified question is open as far as we know:

\begin{Question}\label{degreebound_determinantal}\rm Let $\mathcal{M}$ denote a
square matrix over $k[\xx]=k[x_0,\ldots,x_n]$ whose entries are forms of equal degrees and
such that $f:=\det(\mathcal{M})$ is an irreducible homaloidal polynomial.
Is $\deg(f)\leq n+1$?
\end{Question}
It looks like the question is wide open even if the entries are linear.
In \cite{CRS} a class of examples has been studied in which the answer to the last question
is affirmative.
One main motivation of this work is the hope to shed additional light into this facet of the theory.

 \medskip

 Here is a brief description of the contents of the paper.

 In Section 2 we initially review a  few basic notions of ideal theory to be used throughout. These are mainly related to the syzygies of an ideal and its main associated algebras, stressing some or other notation that may not be universally accepted.
Next there are a some preliminaries about homaloidal polynomials and displays some elementary examples.
A discussion of birationality criteria is annotated for the reader's convenience. Our source for these criteria is generally \cite{AHA}, where the language looks quite encompassing both for the algebraist and the geometer.
Special cases and more particular situations had been exposed before in various other sources (\cite{cremona}, \cite{bir2003}, \cite{CRS}).

Sections 3 and 4 contain the core of the results and  each is subdivided in several subsections.
We focus on determinants of generic matrices of catalectic nature and their degenerations -- we use ``generic" both technically and informally, and in the latter case the emphasis is that the entries are variables, but in the degenerations one allows zeros instead of variables (but no arbitrary linear forms will be entries).
This has the advantage of making the shape of the partial derivatives as close as possible to the submaximal minors.

Section 3 deals with those classical matrices in their generic versions, where entries are always variables.
We start with the completely generic and symmetric generic matrices, as a model to guide us through the other sorts of catalecticants.
Most of what we will say in the fully generic  case  seems to be well-known, still our dealing  recasts what is to be retained from an algebraic angle.
Inspired by the behavior in this case and by the geometric notion of parabolic points, we reinstate a notion (or principle) of parabolism with a more algebraic content solely attached to the Hessian.
For catalecticants the property of parabolism in this sense  turns out to be quite surprising (to us) and even motivated our daring state a conjecture.
Whatever direction the conjecture eventually takes up, the relationship to the ideal theoretic nature of the dual variety is curious and has some measure of puzzlement.

Returning to the homaloidal search, the first consideration is the generic Hankel matrix of arbitrary size.
Thus, let ${\cal H}_m$ stand for the generic $m\times m$ Hankel matrix and let $P\subset R$ denote
its ideal of submaximal minors (i.e., $(m-1)$-minors). We prove that $P$ is the minimal primary
component of the gradient ideal  $J=J(f)\subset R$ of $f:=\det {\cal H}_m$.
Although this assertion sounds naturally guessed, not so much its proof.
We have used  a bouquet of arguments, ranging from multiplicities to initial ideals  to
Pl\"ucker relations (straightening laws) of Hankel maximal minors via  the Gruson--Peskine change of matrix trick.

From this it is but one step to guess that the only remaining associated (necessarily embedded)
prime of $R/J$ is the obvious candidate $Q:=I_{m-2}({\cal H}_m)$ -- defining the singular locus of the determinantal variety $V(P)$.
This guess is equivalent to the expected equality $J:P=Q$.
We chose to include the latter guess in a conjecture of much larger scope to the effect that
$$JP^i:P^{i+1}=I_{m-2-i}({\cal H}_m),$$
for  $0\leq i\leq m-2$.
An almost immediate consequence of this conjectured statement is that $J$ is a minimal reduction
of $P$ with reduction number $m-2$.
This gives quite a spectacular relationship of algebraic content between these two ideals.
Terse as it may look to a geometer, we expect it to pave the road to the next conjecture to the effect
that the gradient ideal $J$ is an ideal of linear type, i.e., that its Rees algebra coincides with its symmetric algebra.
Geometrically, one is saying that the defining equations of the blowup of $\pp^{2(m-1)}$ along the singular locus
of the Hankel determinantal hypersurface can be taken to be of linear nature.

If this conjecture proves to be affirmative then one can conclude, via the criterion stated in Proposition~\ref{polar_and_lineartype}, that the Hankel determinant is not homaloidal.
For this one draws on the result proved in \cite[Theorem 3.3.5]{M_thesis} to the effect that the {\em linear syzygy rank} of
$J$ is $3<2m-2$ for $m\geq 3$.
These various aspects of the nature of the singular ideal give an idea of how much more difficult it is as compared with its nicely behaved minimal component $P$, the latter being linearly presented and thereby defining a birational map onto the image.

So much for the Hankel matrix of arbitrary size.
In the case where $m=3$ we are able to solve all previous conjectures in the affirmative.
The most laborious is the one about the linear type property which, for this value of $m$,
allows to apply some well established criteria in terms of Fitting ideals and Koszul homology.
After this, we move on to other catalecticants of higher leap (step).
The previous results in the Hankel case are not immediately adjustable to these matrices and, in fact, the theory takes a more sinuous route.
First, the ideal of $2$-minors of the $2$-leap $3\times 3$ catalecticant  is not prime, only radical.
One of its minimal primes is generated by certain variables, while the other is the ideal of maximal minors
of an associated catalecticant introduced as an analog of the classical Gruson--Peskine trick.
This already tells us that the minimal part of the gradient ideal is  reduced splitting into two components, thus indicating a measure of difficulty in this case as compared to the Hankel case.
We next deal with the $3$-leap and $2$-leap $4\times 4$ catalecticants.
To derive whether these are homaloidal or not requires quite a bit of juggling in-between theory and computation.
In order to state what is expected in arbitrary dimension and leap, there will be unavoidably some uninspired guessing.
On the bright side, some of the guessing is firmly based on the methods throughout, of which counting the rank of the linear syzygy part of the gradient ideal appeals to some muscle action.
Looking from such an angle, we state the conjecture that only $m$-leap and $(m-1)$-leap $m\times m$ catalecticants will have enough linear syzygy rank.
Unfortunately, we will not get any answers with the computer for $m\geq 5$ as even the computation of the syzygies of the partial derivatives drags along quite a bit.

\smallskip

In Section 4 our focus is on certain degenerations of the generic versions of the previous
section. We avoid calling them specializations since the numerical invariants of the various ideals in consideration
may change as do their properties of interest to this work.
Actually, the only sort of degeneration introduced here is by means of replacing some variables (entries) in strategic
positions by zeros. This idea was originally introduced in \cite[Section 4.1]{CRS} for Hankel matrices.
A reminiscent mention  of homaloidal systems of maximal minors of degenerations of $m\times (m+1)$ Hankel matrices had appeared
in \cite[Section 3]{ST} as a limiting process to solve a problem related to the reciprocal transformation.
The examples show that degenerating can both destroy or give rise to the property of homaloidness.
This phenomenon is yet to be understood from the geometric point of view, but even the underlying algebra is not cleat either.
By and large it is a provocative section which we found appropriate to enlighten ourselves if not the sympathetic reader.

\medskip

The paper includes several results with complete conceptual proofs. On the other hand,  quite a few rely on a mix of conceptual and computational arguments, some of which are followed by various conjectured statements.
It is perhaps fair to say that a facet of this work is to offer an open invitation for others to propose and carry completely conceptual proofs of some of the conjectured statements.  

\medskip

A good deal of the material dealing with Hankel and sub-Hankel matrices was subject of the first author's PhD thesis under the second author's supervision.

We wish to thank S. H. Hassanzadeh for his help and suggestions.

\section{Preliminaries}

\subsection{Tools from ideal theory}

Let $R$ be a Noetherian ring and  let $I\subset R$ be an
ideal.

 Let ${\cal S}_R(I)\surjects {\cal R}_R(I)$ denote the structural
graded $R$-algebra homomorphism from the symmetric algebra of $I$ to
its  {\sc Rees algebra} - the latter is the graded $R$-algebra that defines the blowup
along the subscheme corresponding to the ideal $I$ (see \cite[\S
5.2] {Eisenbook}). We say that $I$ is of {\sc linear type\/} if this
map is injective.

An ideal $I\subset R$ of linear type satisfies  the {\sc
Artin--Nagata condition $G_{\infty}$\/}  (see \cite{AN}) which states
that the minimal number of generators of $I$ locally at any prime
$p\in \spec{R/I}$ is at most the codimension of $p$. It is known that this condition is
equivalent to a condition in terms of a free presentation
$$R^m\stackrel{\varphi}{\lar} R^{n+1} \lar I \lar 0 $$
of $I$ and the Fitting ideals of $I$, namely:
\begin{equation}\label{f1}
\hht (I_t(\varphi))\geq {\rm rank}(\varphi)-t+2,\quad \mbox{\rm
for}\quad 1\leq t\leq {\rm rank}(\varphi),
\end{equation}
where $I_t(\varphi)$ denotes the ideal generated by the $t\times
t$ minors of a representative matrix of $\varphi$ and $\hht$ designates the {\sc height}
of an ideal (see, e.g., \cite[\S 1.3]{Wolmbook}).
Because of its formulation in terms of Fitting ideals, the condition has
been dubbed as {\sc property} $(F_1)$.

 Suppose that $R$ is a standard graded over a field $k$ and $I$ is
 generated by forms of a given degree $s$. In this
case, $I$ is more precisely given by means of a free graded
presentation
$$R(-(s+1))^{\ell}\oplus\sum_{j\geq 2} R(-(s+j)) \stackrel{\varphi}{\lar} R(-s)^{n+1}\rar
I\rar 0$$
 for suitable shifts $-(s+j)$ and rank $\ell\geq 0$.
 Of much interest in this work is the value of $\ell$, so let us state in which form.
We call the image of $R(-(s+1))^{\ell}$ by $\varphi$ the {\sc linear
part\/} of $\varphi$ -- often denoted $\varphi_1$. One says that the rank of $\varphi_1$
is the {\sc linear rank} of $\varphi$ (or of $I$ for that matter) and that $\varphi$ has {\sc maximal linear rank\/}
provided its linear rank is $n$ ($={\rm
rank}(\varphi)$). Clearly, the latter condition is trivially
satisfied if the $\varphi_1=\varphi$, in which case $I$ is said to have
{\sc linear presentation\/} (or is {\sc linearly presented\/}).

A quite notable fact, whenever $R=k[\xx]=k[x_0,\ldots,x_n]$ is a standard graded polynomial ring over $k$,
is the case where $I$ happens to be of linear type and generated by $r+1$ forms of the same degree.
Then $I$ has maximal analytic spread and hence ($k$ infinite) is
generated by analytically independent forms (with respect to the irrelevant maximal ideal) of the same degree.
 Then these forms are algebraically independent elements over $k$, hence define a dominant rational map $\p^n\dasharrow \p^n$.
 This will be a cornerstone of many an argument to follow.

\subsection{Birationality criterion in terms of ideals}

Let $f\in k[\xx]=k[x_0,\ldots,x_n]$ be a squarefree homogeneous
polynomial of degree $d\geq 2$. Let
$$I=\left(\frac{\partial f}{\partial x_0},\ldots,\frac{\partial f}{\partial x_n}\right)
\subset k[\xx],$$ the so-called {\sc gradient\/} ideal of $f$ -- we
refrain from using the terminology {\sc Jacobian\/} ideal since the
latter usually refers to the residue ideal in the ring $k[\xx]/(f)$.

The partial derivatives of $f$ can be looked upon as the coordinates
of a rational map ${\cal P}_f\colon\pp^n\dasharrow \pp^n$.  This map
is called the {\sc polar map\/} defined by (or of) $f$. We note that the
image of this map is the projective subvariety on the target whose
homogeneous coordinate ring is given by the $k$-subalgebra
$k[{\partial f}/{\partial x_0},\ldots,{\partial f}/{\partial
x_n}]\subset k[\xx]$ (whose grading, for that purpose, is renormalized so as to
have its generators of degree one). While this algebra describes the
image of the map it falls short of giving the complete picture of
the map itself. A thorough analysis of the polar map has been
undertaken in \cite{CRS}. Here we focus on the algebraic properties
of this map and algebraically structured examples.

Note that the ideal $I$ has codimension at least two since $f$ is
assumed to be squarefree. In particular, the partial derivatives
give the unique representative of codimension $\geq 2$ of the polar
map (cf. \cite[Proposition 1.2]{bir2003}). Thus, it makes sense to
call the singular scheme of $f$ the (uniquely defined) {\sc base
scheme\/} of the polar map. A complete understanding of the nature
of the polar map hinges on describing its base scheme.

If ${\cal P}_f$ is birational one says that $f$ is {\sc homaloidal}.
One may expect that the finer properties of ${\cal P}_f$ are
embodied in a rich interplay between the geometric side of the map
and the properties of the gradient ideal $I$.
In this vein, we will refrain from mechanically pass to
the projective scheme defined by the gradient ideal $I$ since
knowing whether $I$ is $(\xx)$-saturated may have direct bearing to
the properties of the polar map.

\medskip

A general characteristic-free birationality criterion has been established in \cite{AHA}
which depends on a unique numerical invariant.
This invariant can be viewed as a replacement for the field (topological) degree of the given map and, for many purposes,
it is more flexible and computationally effective.
For polar maps such that the corresponding linear system generates an ideal of linear type,
the following simplified criterion holds in all characteristics.

\begin{Proposition}\label{polar_and_lineartype} {\rm \cite[Theorem 3.2 and Proposition 3.4]{AHA}}
Let $f\in R:=k[\xx]=k[x_0,\ldots,x_n]$ denote a square-free homogeneous polynomial of degree
$d\geq 3$ such that the image of the corresponding polar map has dimension $n$.
Consider the following conditions statements:
\begin{enumerate}
\item[{\rm (a)}] The syzygy matrix
of the gradient ideal of $f$ has maximal linear rank.
\item[{\rm (b)}] $f$ is homaloidal.
\end{enumerate}
Then {\rm (a)} implies {\rm (b)}. If, moreover, the gradient ideal is of linear type then {\rm (a)}
and {\rm (b)} are equivalent.
\end{Proposition}

As a reminder, an ideal of linear type generated in fixed degree is generated by algebraically independent forms, hence the rational map defined by these forms is dominant (onto $\pp^n$).
As the converse is false, this is a source of difficulty to go from the geometric to the algebraic aspect.

Additional impact to the algebraic side is a certain polynomial akin to the classical {\em principal curves}
in plane Cremona map theory, being related to the  so-called {\em Jacobian curve} of a homaloidal net.
Its terminology has been formally introduced in
\cite[Propositions 1.2 and 1.3]{Zaron} for the
use in the theory of symbolic powers, while
earlier appearances are in \cite[Proof of Theorem 3.1]{CremonaMexico} (see also \cite{CostaSimis})
for the case of maps defined by monomials.
The object itself is part of the algebraic nature of a Cremona map, in that it ties the map and its inverse
by means of an essentially unique polynomial.
One interest in the present work is as to how this polynomial in the case of a homaloidal
polynomial is further related to the polynomial itself.

We briefly recall the definition.
Let $k$ denote an arbitrary infinite field  -- further assumed to be algebraically closed in a geometric discussion.
Recall that, quite generally, a rational map $\mathfrak{F}:\pp^n\dasharrow \pp^n$ is defined by $n+1$ forms $\mathbf{f}=\{f_0,\ldots, f_n\}
\subset R:=k[\xx]=k[x_0,\ldots,x_n]$ of the same degree $d\geq 1$, not all null.
We often write $\mathfrak{F}=(f_0:\cdots :f_n)$ to underscore the projective setup
and assume that $\gcd\{f_0,\cdots ,f_n\}=1$ (in the geometric terminology, the linear system defining $\mathfrak{F}$ ``has no fixed part''),
in which case we call $d$ the {\em degree} of $\mathfrak{F}$.

If $\mathfrak{F}$ is a Cremona map of $\pp^n$ then there is a rational map
 $\mathfrak{G}:\pp^{n}\dasharrow \pp^{n}$ based on a linear system spanned by forms $\mathbf{g}=\{g_0,\ldots, g_n\}
\subset R$ of same degree  satisfying the
relation
\begin{equation}\label{birational_rule}
(\mathbf{g}_0(\mathbf{f}):\cdots :\mathbf{g}_n(\mathbf{f}))\equiv (x_0:\cdots :x_n).
\end{equation}
The congruence translates into the existence of a uniquely defined form $D\in R$ such that, using a short vector notation,
$\mathbf{g}(\mathbf{f})=D\cdot (\xx)$.
In \cite{Zaron} $D$ has been dubbed the {\em inversion factor} of the map $\mathcal{F}$ or, more precisely, its
{\em source inversion factor}.
Its degree is  $\deg (\mathfrak{F})\deg(\mathfrak{G})-1$.

For convenience, we quote the following basic result about the inversion factor:

\begin{Proposition}\label{jac_vs_factor}{\rm (\cite[Proposition 1.3]{Zaron}, char$(k)=0)$}
Let $\mathfrak{F}$ denote a Cremona map of $\pp^{n}$ defined by forms $\ff:\{f_0,\ldots,f_n\}\subset R$
without fixed part and let $\Theta(\ff)$ denote the Jacobian matrix of $\ff$.
Then $\det \Theta(\ff)$ divides a power of the source inversion factor $E$ of $\mathfrak{F}$.
In particular, if $\det \Theta(\ff)$ is reduced then it divides $E$.
\end{Proposition}

\section{Matrices. I: generic catalecticants and symmetric}

Since we are mainly interested in the search of irreducible homaloidal polynomials, is is
natural to start looking about in determinants of well-structured square matrices.
All matrices in this section, except for the symmetric ones, are of the following type.

Let $m\geq 2$ and $1\leq r\leq m$ be given integers.
Let $R=k[x_0,\ldots,x_n]$ be a polynomial ring with $n+1=(m-1)(r+1)+1$ variables.
The  $r$-{\em leap} $m\times m$  {\em generic catalecticant} is the matrix
\begin{equation}\label{gen_catalectic}
\mathcal{C}_{m,r}=
\begin{pmatrix}
X_0 & X_1 & X_2 & \hdots & X_{m-1}\\
X_{r} & X_{r+1} & X_{r+2} & \hdots & X_{m+r-1}\\
X_{2r} & X_{2r+1} & X_{2r+2} & \hdots & X_{m+2r-1}\\
\vdots & \vdots & \vdots & \ddots & \vdots\\
X_{(m-1)r} & X_{(m-1)r+1} & X_{(m-1)r+2} & \hdots & X_{(m-1)r+(m-1)}
\end{pmatrix}
\end{equation}
Note that the corresponding
determinant will have low degree ($=n$) as compared to the dimension of the
ring and still involve all variables.

The extreme values $r=1$ and $r=m$ yield, respectively, the ordinary Hankel matrix and the generic matrix.
An important result proved in \cite{Eisenbud2} is that the generic Hankel matrix of arbitrary size $m\times n$ is $1$-generic.
Using this result, it has been proved in \cite{Zaron2} that all generic catalecticants  of arbitrary size are $1$-generic.
The advantage of this notion is that it implies, in particular, that the  determinant of a square such a matrix is irreducible.

Although the fully generic matrix is an extreme case of a catalecticant, not so the generic symmetric matrix.
However, the two have similar behavior regarding  the subject of this work, and hence will be treated together.

\subsection{Generic and generic symmetric matrices}\label{symmetric}

The main part of the following result is classically known in algebraic geometry.
We restate it by emphasizing the algebraic side, along with a tiny addition.

Recall that a homogeneous polynomial is called {\em totally Hessian}
if up to a nonzero element of $k$ its Hessian determinat is a power of it.
Thus, if $f\in k[x_0,\ldots,x_n]$ has degree $d\geq 1$, being totally Hessian means that
$$H(f)=cf^{\frac{(d-2)(n+1)}{d}},$$
where $H(f)$ denotes the Hessian of $f$ and $c\in k\setminus \{0\}$.
The totally Hessian situation is to be considered as an extremal condition, whereby one might expect that
more commonly a certain power of $f$ divide the Hessian determinant.
For reasons that will be explained soon, the maximal exponent for which a power of $f$ may happen to divide $H(f)$ is
said to be the {\em expected multiplicity} if its value is $n-1-v(f)$, where $n$ is the dimension of the ambient projective space and
$v(f)=\dim V(f)^*$ (the dimension of the dual variety to $V(f)$).
For further terminology we refer to \cite[Section 2.1]{CRS}, which contains a discussion of these notions referring back to B. Segre.

\begin{Proposition}\label{genericdet}
Let $\xx$ denote a generic or a generic symmetric $m\times m$  matrix over the field $k$
and let $f=\det\xx$.
Then:
\begin{enumerate}
\item[{\rm (a)}] $f$ is a  homaloidal polynomial and the polar map of $f$ is an involution up to a projective transformation.
\item[{\rm (b)}] The source inversion factor of the polar map of $f$ coincides with the $(m-2)$th power of $f$.
\item[{\rm (c)}] $f$ is totally Hessian$\,${\rm ;} in particular,
the expected exponent of $f$ as a factor of $H(f)$ has the expected value if and only if $m=3$.
\end{enumerate}
\end{Proposition}
\demo We first remark that in the generic case the partial derivatives of $f$ are the signed
cofactors of the matrix $\xx$, which is immediate from the data.
In the case of the generic symmetric matrix each $(m-1)$-minor appears twice, hence the partial derivatives relative to the variables off the main diagonal will be the corresponding cofactor multiplied by $2$.

(a) The ideal of $k[\xx]$ generated
by the cofactors is of linear type and linearly presented (cf.
\cite{Huneke} for the generic case and \cite{Kotsev} for the
generic symmetric case (cf. also \cite{Wolmbook}).
Therefore, the assertion about homaloidness is a consequence of Proposition~\ref{polar_and_lineartype}.

Since $f$ read in dual variables is the dual hypersurface to
the variety defined by the cofactors ($(m-1)$-minors) then the inverse to the polar map of $f$ is of the same kind
up to a linear change of coordinates.
(Note that what is a true involution is the map defined by the cofactors, hence in the fully generic case the polar map is a true involution, while in the symmetric case it will be so up to suitable coefficient scrambling.)

(b) We argue in the generic case, the symmetric case requiring small adjustments.
Let $\underline{\Delta}$ denote the ordered list of the signed $(m-1)$-minors (i.e., the cofactors of the matrix).
As noted above, it defines an involution, hence one has by (\ref{birational_rule}):
\begin{equation}\label{involution_eq}
\Delta_{x_i}(\underline{\Delta})=x_iD,
\end{equation}
for every $0\leq i\leq m^2$, where $\Delta_{x_i}$ denotes the (signed) cofactor of $x_i$ and $D$ is the inversion factor.
Using Laplace to compute $f$ in the form
$f=\sum_{j=0}^{m-1}x_j\Delta_{x_j},$ and applying (\ref{involution_eq}), one obtains
one obtains
\begin{equation}\label{basic_evaluation}
f(\underline{\Delta})=\sum_{j=0}^{m-1}\Delta_{x_j}\Delta_{x_j}(\underline{\Delta})=\sum_{j=0}^{m-1}\Delta_{x_j}\,x_jD= D\,\sum_{j=0}^{m-1}x_j\Delta_{x_j}=D\,f.
\end{equation}
Then $D=f^{m-2}$ provided $f(\underline{\Delta})=f^{m-1}$.
But the latter follows from Cauchy's formula for the cofactors which asserts that  the determinant of the matrix of cofactors $C(\xx)$ (the so-called adjugate of $f=\det (\xx)$) is equal to $f^{m-1}$.
Since  $(\det (\xx)) (\underline{\Delta}) = \det( \xx (\underline{\Delta}) )=\det (C(\xx)^t)=\det (C(\xx))$, we are through (we thank Z. Ramos for pointing to us this  passage).

(c)
For the assertion about being totally Hessian,  we have to show that $H(f)=f^{m(m-2)}$.
By Proposition~\ref{jac_vs_factor}, $H(f)$  divides a power of the inversion factor $D=f^{m-2}$, hence is itself a power of $f$ as $f$ is irreducible. A degree count gives that $H(f)=f^{m(m-2)}$.

\smallskip

As for the subsumed assertion, in the present situation, $n=m^2-1$.
As remarked in the proof of (a), $V(f)^*\subset \pp^{m^2-1}$ is the subvariety defined by the $(m-1)$-minors of the matrix (read in dual variables). But the latter has codimension $4$.
therefore, we get
$$n-1-v(f)=m^2-2-(m^2-1-4)=3.$$
On the other hand, since $f$ is totally Hessian, the true exponent of $f$ as a factor of $H(f)$ is $m(m-2)$. Clearly, $m(m-2)=3$ if and only if $m=3$.
This proves the subsumed assertion.
\qed

\begin{Remark}\rm
(1) A shorter elementary proof of part (a), more in the spirit of the linear algebra used in the proof of (b), has been communicated to us by F. Russo.

(2) There has been interest in considering polynomials $f$ which are factors of their Hessian determinant with multiplicity higher that the expected one (see \cite[Section 2.1]{CRS} for appropriate references.)

(3) It would be curious to decide if there exist strict linear specializations of the generic matrix -- that is to say, $n\times n$ matrices of linear forms in less than $m^2$ variables
whose  $(m-1)$-minors still generate a codimension $4$
Gorenstein ideal -- such that the dual variety to the $(m-1)$-minors is a homaloidal hypersurface.
Though some of these will have a homaloidal determinant, 
the dual variety may even fail to be a hypersurface (see Remark~\ref{dual_variety}).

A similar question can be asked in the case of specializations of the generic symmetric matrix.
\end{Remark}

\subsection{The Hessian principle: parabolism}

The idea of this part is to relate more closely the two hypersurfaces $V(f)$ and $V(H(f))$, where $f\in R=k[x_0,\ldots, x_n]$ is the determinant  of certain specializations of the generic and the generic symmetric matrices.

Assuming that $f$ is irreducible, we follow the terminology established in \cite[Section 2.1]{CRS} taken from the references mentioned there.
Accordingly, to say that the generic point of $V(f)$ is  $h$-parabolic, for some integer $h> 0$, translates into $f$ being a factor of $H(f)$ with multiplicity at least $h$; then, necessarily, $h=n-1-\dim V(f)^*$.
Note that one cannot decide the non-vanishing of $H(f)$ via this definition, since if $H(f)=\{0\}$ then certainly $f$ is a factor of $H(f)$ with any multiplicity (cf. \cite[Corollary 4.4]{Zak} and the comments thereafter for this puzzling issue).

Drawing upon this line of ideas, we will say that a homogeneous polynomial $f\in R$ is {\em parabolic} if $f$ is a factor of $H(f)$ with bounded multiplicity (not necessarily the expected multiplicity).
This means, in particular, that $H(f)\neq 0$ and that the (true) multiplicity of $f$ as a factor of $H(f)$ is an upper bound for $n-1-\dim V(f)^*$.
(Calling parabolic even if the multiplicity is not the expected one seems like a good idea as the two properties can ba handled separately.)

\begin{Conjecture}\rm
Let $f$ denote the determinant of the $r$-leap $m\times m$ catalecticant  (\ref{gen_catalectic}), with $1\leq r\leq m$.
Then $f$ is parabolic and if $r\leq m-1$ then it has the expected multiplicity as a factor of $H(f)$.
\end{Conjecture}

As mentioned before, $f$ is irreducible since the generic catalecticant is $1$-generic (\cite[Proposition 2.1]{Zaron2}).
Note that $r\leq m-1$ excludes the case of the fully generic matrix (though it could be included if $m=3$).

Let us consider the case $r=m-1$, which is the case where one expects the Hessian determinant $H(f)$ to admit $f$ as factor of highest multiplicity among the class of $r$-leap catalecticants with $r\leq m-1$.

As a slight, but beautiful, evidence to the conjecture, we can entirely describe the prime factorization of $H(f)$ for $m=3,4$:

\begin{Example}\label{parabolic_exs}\rm
$\bullet$ ($m=3$) One is looking at the matrix
$$\mathcal{C}_{3,2}=\left(
  \begin{array}{ccc}
    x_0 & x_1 & x_2 \\
    x_2 & x_3 & x_4 \\
    x_4 & x_5 & x_6
  \end{array}
\right).
$$
Here one has $H(f)=f\cdot g$, where $g$ is (up to a projective change of coordinates) the equation of the dual surface to the twisted cubic in the variables $x_0,x_2,x_4,x_6$.
Note that the defining equations of the latter are the $2\times$2 minors of the first and third columns of $\mathcal{C}_{3,2}$.
The projective change of coordinates is totally trivial, obtained by $x_0\mapsto 3x_0, x_6\mapsto 3x_6$ and fixing the remaining variables.
Observe that the degrees match as the equation of the dual variety has degree $4$ and $H(f)$ has degree $7$.

\smallskip

The dual to the hypersurface $V(f)$ is an arithmetically Cohen--Macaulay variety of codimension $2$,  a
determinantal scheme defined by the maximal minors of a linear $4\times 3$ matrix.
This scheme can be obtained by the method of \cite[Proposition 1.1 (ii)]{CRS} as the projection
from $\pp^8={\rm Proj}(k[z_0,\ldots,z_8])$ to $\pp^6$ with projecting center cut
by the linear forms $z_0,z_1,z_2-z_3,z_4,z_5-z_6,z_7,z_8$ (note that these define the specialization from the generic case to
the present catalecticant case).

The expected multiplicity of $f$ is $n-1-\dim V(f)^*=5-4=1$, which coincides with the effective multiplicity above.

\medskip

$\bullet$ ($m=4$) One is looking at the matrix
$$\mathcal{C}_{4,3}=\left(
  \begin{array}{cccc}
    x_0 & x_1 & x_2  & x_3\\
    x_3 & x_4 & x_5  & x_6 \\
    x_6 & x_7 & x_8 & x_9 \\
    x_9 & x_{10} & x_{11} & x_{12}
  \end{array}
\right).
$$
Here one has $H(f)=f^5\cdot (\det \mathcal{H}_3)^2$, where $\mathcal{H}_3$ is the Hankel matrix
$$\mathcal{H}_{3}=\left(
  \begin{array}{ccc}
    x_0 & x_3 & x_6 \\
    x_3 & x_6 & x_9 \\
    x_6 & x_9 & x_{12}
  \end{array}
\right).
$$
Note that the $2\times 4$ matrix corresponding to $\mathcal{H}_3$ by the classical principle (\cite[Lemme 2.3]{GP}) is the submatrix of $\mathcal{C}_{4,3}$ with first and fourth columns.
Observe that the degrees match.
Note, however, that this time around the dual variety of the normal quartic in variables $x_0,x_3,x_6,x_9,x_{12}$ is not coming into the picture --  it has the right matching degree (=$6$) but is irreducible.

\smallskip

The dual to the hypersurface $V(f)$ has codimension $6$.
The expected multiplicity of $f$ is $n-1-\dim V(f)^*=11-6=5$, which coincides with the effective multiplicity above.
\end{Example}

Let us now consider examples at the other extreme of the value of $r$, namely, for  Hankel matrices.

\begin{Example}\label{parabolic_hankel}\rm
$\bullet$ ($m=3$) One is looking at the matrix
$$\mathcal{H}_{3}=\left(
  \begin{array}{ccc}
    x_0 & x_1 & x_2 \\
    x_1 & x_2 & x_3 \\
    x_2 & x_3 & x_4
  \end{array}
\right).
$$
The prime factorization of the Hessian determinant is $H(f)=f\cdot \partial$, where $\partial$ denotes the partial derivative  $\partial f/\partial x_2$ up to a sign and a trivial projective change of coordinates.

The dual variety $V(f)^*$ is a well-known codimension $2$ subvariety of $(\pp^4)^*$ defined by $7$ cubics, obtained either by suitable projection of the $2$-Veronese or else as image of the rational map defined by a linear system of quadrics of dimension $4$ in $3$ variables.
Then, the expected multiplicity is $n-1-\dim V(f)^*=3-2=1$, clearly coinciding with the effective multiplicity above.

\smallskip

$\bullet$ ($m=4$) One is looking at the matrix
$$\mathcal{H}_{4}=\left(
  \begin{array}{cccc}
    x_0 & x_1 & x_2  & x_3\\
    x_1 & x_2 & x_3  & x_4 \\
    x_2 & x_3 & x_4 & x_5 \\
    x_3 & x_4 & x_5 & x_6
  \end{array}
\right).
$$
One finds that the effective multiplicity of $f$ is $2$. A reasonable bet is that the complementary factor is the square of $\partial f/\partial x_3$
up to an appropriate projective change of coordinates, in which case the prime factorization would be
$H(f)=f^2\cdot (\partial f/\partial x_3)^2$ up to coordinate change.
Unfortunately, we have not been able to decide if this is the case or if this factor of degree $6$ is irreducible.

The dual variety of $V(f)$ is a codimension $3$ subvariety of $(\pp^6)^*$ defined by  quartics.
As it comes out, the expected multiplicity $5-3=2$ coincides with the effective multiplicity.
\end{Example}

Finally, we consider one more case, where $r$ has an intermediate value so to say.

\begin{Example}\label{parabolic_cat2}\rm
$$\mathcal{C}_{4,2}=\left(
  \begin{array}{cccc}
    x_0 & x_1 & x_2  & x_3\\
    x_2 & x_3 & x_4  & x_5 \\
    x_4 & x_5 & x_6 & x_7 \\
    x_6 & x_7 & x_8 & x_9
  \end{array}
\right).
$$
The effective multiplicity of $f$ is $2$.
The complementary factor (of degree $12$) is pretty unreachable for inspection, except that it belongs to the gradient ideal. It is plausible that, up to a coordinate change, it is a product of partial derivatives of $f$.

The dual $V(f)^*$ is a codimension $3$ arithmetically Cohen--Macaulay subvariety of $(\pp^9)^*$ with linear resolution.
By a well-known result (\cite[HuMi]), the expected multiplicity is $8-6=2$, which is the value of the effective one.

\end{Example}
\begin{Remark}\rm
All the required computer calculation has been performed with {\em Macaulay} (\cite{Macaulay}).
Cases beyond these will most probably get stalled due to the computation of determinants and dual varieties of hypersurfaces in large ambient spaces. A breach for the theory to come in is a specialization procedure from the fully generic case. Care must be taken, since although the specialization of the Hessian of the generic matrix contains as a submatrix the Hessian matrix of the specialized generic matrix (i.e., the catalecticant), the determinant of the latter won't be a factor of the former.

As for the theoretic tool to get to the expected multiplicity at the moment we have no insight as far as algebraic methods are concerned. It is possible that some of the methods employed in \cite{CRS} might be used.
\end{Remark}

\subsection{Hankel matrices}

Hankel matrices constitute the lower extreme of a catalecticant in terms of its leap.
For this reason and also because they are the only symmetric catalecticants we deal with them first.

 In general, an obstruction for the
cofactor technique, which has shown so successful in the fully generic and generic symmetric cases, is that the minimal number of generators of the ideal generated by the cofactors has to equal the number of
variables retained in the specialized matrix.
This immediately rules out the cofactor technique for $r$-leap catalecticant square matrices with $r\leq m-1$.  To
rephrase it, the problem is that the gradient of the generic
determinant does not specialize to the gradient of the specialized
determinant.

The generic Hankel matrix of size $m\times m$ is the symmetric matrix
$$\mathcal{H}_m:=
 \left(
\begin{matrix}
x_0&x_1&\ldots &x_{m-1}\\
x_1&x_2&\ldots &x_{m-2}\\
\vdots &\vdots &\ldots &\vdots \\
x_{m-1}&x_{m-2}&\ldots &x_{2m-2}\\
\end{matrix}
\right)
$$
Let $J \subset R=k[x_0,\ldots, x_{2m-2}]$ denote the gradient ideal of $\det \mathcal{H}_m$ and let $P:=I_{m-1}(\mathcal{H}_m)\subset R$ stand for the ideal of $(m-1)$-minors of the matrix.

Our main goal in this subsection is to shows several properties relating these two ideals.
Throughout the multiplicity (degree) of a graded residue ring $R/I$ is denoted $e(R/I)$.

\begin{Lemma}\label{degree} Set $P:=I_{m-1}(\mathcal{H}_m)$. Then $P$ is a prime ideal and the multiplicity of $R/P$ is
$$e(R/P)=\frac{1}{3!} (m-1)m(m+1)$$
\end{Lemma}
\demo By \cite[Proposition 4.3]{Eisenbud2}, the ideal $P$ is prime and has maximal possible height.
On the other hand, by the same principle of \cite[Lemme 2.3]{GP}, $P$ is the ideal of the maximal minors of the $(m-1)\times (m+1)$ Hankel matrix
\begin{equation}\label{GP_trick}
 \left(
\begin{matrix}
x_0&x_1&\ldots &x_m\\
x_1&x_2&\ldots &x_{m-1}\\
\vdots &\vdots &\ldots &\vdots \\
x_{m-2}&x_{m-1}&\ldots &x_{2m-2}\\
\end{matrix}
\right)
\end{equation}
Therefore, this ideal has height $m-(m-2)+1=3$.
It follows that the Eagon--Northcott complex resolves $R/P$; since this complex is well-known to be a pure $(m-1)$-resolution
one can apply the formula of \cite{HuMi} to get the desired expression.
\qed

\begin{Proposition}\label{codimension} The gradient ideal $J$ of $\det \mathcal{H}_m$ is a codimension $3$ ideal contained in $P$.
\end{Proposition}
\demo There are several ways to observe the containment $J\subset P$.
We quote a more general fact which yields this result:

\begin{Lemma}\label{partials_vs_cofactors}{\rm \cite{Golb}, \cite{M_thesis}}
Let $\mathcal{M}$ denote a square matrix over $R=k[z_0,\ldots,z_s]$ satisfying the following requirements:
\begin{itemize}
\item Every entry of $M$ is either $0$ or $z_i$, for some $i=0,\ldots, s$
\item Any variable $z_i$ or $0$ appears at most once on every row or column.
\end{itemize}
Let $f:=\det(\mathcal{M})\in R$. Then, for each $i=0,\ldots, m$, the partial derivative $f_i$ of $f$ with respect to $z_i$ is the sum
of the {\rm (}signed{\rm )} cofactors of the entry $z_i$ in all of its appearances on $\mathcal{M}$.
\end{Lemma}

\medskip

We conclude the proof of the Lemma by showing that $J$ has codimension at least $3$.
For this we consider the initial ideal of $J$ in the reverse lexicographic order.
Using Lemma~\ref{partials_vs_cofactors} with $\mathcal{M}=\mathcal{H}_m$ and $R=k[x_0,\ldots,x_{2m-2}]$,
direct inspection shows that ${\rm in}(f_0)=x_m^{m-1}, \,{\rm in}(f_{2m-2})=x_{m-2}^{m-1},
\, {\rm in}(f_{m-1})=mx_{m-1}^{m-1}$. Clearly  then ${\rm in}(J)$ has codimension at least $3$.
Therefore, $J$ has codimension at least $3$ as well.
\qed

\begin{Lemma}\label{multiplicity_inequality}
Consider the subideal
$J'=({\rm in}(f_0),\ldots, {\rm in}(f_{2m-2})\subset R$ of  the initial ideal of $J$ in
he reverse lexicographic order. Then $e(R/J')<2 e(R/P)$.
\end{Lemma}
\demo Let us compute $e(R/J')$ in a direct way.
Namely, using Lemma~\ref{partials_vs_cofactors} once more with $\mathcal{M}=\mathcal{H}_m$ and $R=k[x_0,\ldots,x_{2m-2}]$,
one sees that $J'$ is generated by the set of monomials
$$G:=\{x_{m-2}^{m-1},\,x_{m-2}^{m-2}x_{m-1},\,\ldots,\,x_{m-1}^{m-1},\,x_{m-1}^{m-2}x_m,\,\ldots,\, x_{m}^{m-1}\}.$$
Writing $J''$ for the ideal of $S:=k[x_{m-2},x_{m-1},x_m]$ generated by these monomials one has $J'=J''R$, hence
$e(R/J')=e(S/J'')$.
Since $S/J''$ is a graded of finite length, $e(S/J'')=\dim_k(S/J'')$.
We claim that
\begin{equation}\label{mult_of_pre_initial}
\dim_k(S/J'')=\frac{1}{6} (m-1)m(2m-1).
\end{equation}
To compute the latter, we argue that a basis thereof is formed by all monomials in $S$ of degrees $\leq m-2$ and in addition
the monomials in degrees $m-1, m,\ldots, 2(m-2)=2m-4$ that are not multiples of the generators of $J''$.
The first compound is given by the well-known number
$$\sum_{0\leq d\leq m-2}{{d+2}\choose {2}}={{m+1}\choose {3}}.$$
For the second part, for every $d\in \{m-2, m-1,m,\ldots, 2m-4\}$, let ${\rm M}_d$ denote the set of monomials of degree $d$ that
are divisible by both $x_{m-2}$ and $x_m$. The cardinality of the complementary
set ${\rm M}_d\setminus (x_{m-1}){\rm M}_{d-1}$
can be seen to be
$${{2m-d-2}\choose {2}}.$$

Adding up for $d=m-1,m,\ldots, 2m-4$ gives
$$\sum_{2m-4\geq d\geq m-1} {{2m-d-2}\choose {2}}=\sum_{2\leq e\leq m-1}{{e}\choose {2}}={{m}\choose {3}}.$$

Summing up we have
$$\dim_k(S/J'')={{m+1}\choose {3}}+{{m}\choose {3}}=\frac{1}{6} (m-1)m(2m-1),$$
as was to be shown.

\medskip

An alternative argument yielding this formula goes as follows:

For lighter reading set $x=x_{m-2}, y=x_{m-1}, z=x_m $. We partition  the basis of $(S/J'')$ into disjoint sets according to the powers of $y$. Let then $A_i$ denote the set of all basis elements with $y$-degree  $i-1$, for $1\leq i \leq m-1.$  In other words, $A_i$ denotes the set of all monomials $x^\alpha ÿy^\beta z^\gamma$ such that $\beta=i-1$, $\alpha<m-i$ and $\gamma< m-i$. Indeed, if to the contrary  $\alpha \geq m-i$ then the summation of the $y$-degree and the $x$-degree is $\geq i-1+m-i=m-1$, so the monomial would belong to  $G.$ The same argument applies to show that, for any $i$, the highest possible value of  $\gamma$ in $A_i$ is $m-i-1.$ By inspecting the nature of the monomials in $G$ it ensues that $\alpha$ and $\gamma$ run independently. It therefore follows that

{\small
\begin{eqnarray}\nonumber
A_i&=&\{ y^{i-1}, y^{i-1}z, \ldots, y^{i-1}z^{m-i-1}, y^{i-1}x, y^{i-1}x z,\ldots , y^{i-1}x z^{m-i-1}, \ldots , y^{i-1}x^{m-i-1},\\
\nonumber
&&y^{i-1}x^{m-i-1}z, \ldots , y^{i-1}x^{m-i-1}z^{m-i-1}\},
\end{eqnarray}
where $1\leq i \leq m-1.$
}

It follows that
\begin{equation*}
\dim_k(S/J'')= \sum^{m-1}_{i=1}|A_i|=\sum^{m-1}_{i=1}(m-i)^2=
\frac{1}{6} (m-1)m(2m-1).
\end{equation*}

\medskip

As a consequence of (\ref{mult_of_pre_initial}) one has
$$e(R/J')= \frac{1}{6} (m-1)m(2m-1)< \frac{1}{3} (m-1)m(m+1)=2 e(R/P),$$
as required.
\qed

\medskip

Next we introduce some ideas from the theory of Pl\"ucker relations.
For this we recast the trick of passing to the matrix (\ref{GP_trick})
whose maximal minors are our $(m-1)$-minors in a suitable disposition.
This allows to express some formulas related to the Pl\"ucker relations of these maximal minors.

First, the set of maximal minors of the matrix (\ref{GP_trick}) are partially ordered
in the usual way, using a well-known notation:

\begin{equation}\label{lex ord minor}
[i_1, \ldots , i_{n-1}]\leq [j_1,\ldots , j_{n-1}]   \Leftrightarrow
  i_1\leq j_1, \ldots , i_{n-1}\leq j_{n-1}.
\end{equation}
As an illustration, in the case of $m=4$, we get

\medskip

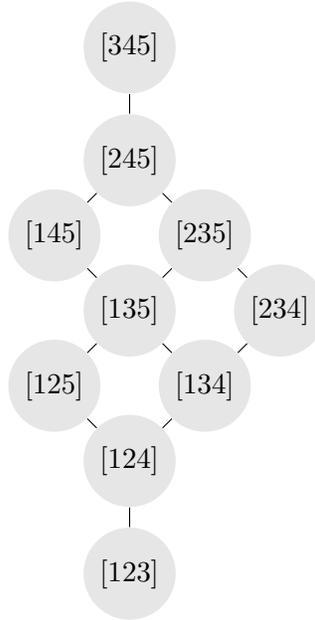
\hspace{2.2in}
\begin{tikzpicture}
  [scale=1.0,auto=center,every node/.style={circle,fill=grey}]
  \node (n1) at (2,0.5) {[123]};
 \node (n2) at (2,2)  {[124]};
  \node (n3) at (1,3)  {[125]};
  \node (n4) at (3,3) {[134]};
  \node (n5) at (2,4)  {[135]};
  \node (n6) at (4,4)  {[234]};
  \node (n7) at (1,5)  {[145]};
  \node (n8) at (3,5)  {[235]};
  \node (n9) at (2,6)  {[245]};
  \node (n10) at (2,7.5)  {[345]};
  \foreach \from/\to in {n1/n2,n2/n3,n2/n4,n3/n5,n4/n5,n4/n6,n5/n7,n5/n8,n6/n8,n7/n9,n8/n9,n9/n10}
    \draw (\from) -- (\to);

\end{tikzpicture}

\medskip

In this diagram the minors horizontally aligned are incomparable.
One important feature of the Pl\"ucker relations is that they can be used to get a straightening law for the minors, whereby a product of incomparable ones is a combination of products of minors each having a factor which is less than the starting minors.
For the contents on this topic in general we refer to  \cite[chapter 4]{BV} and for the present case of Hankel matrices our reference is \cite{Conca}.

From the original square Hankel matrix $\mathcal{H}_m$, by \cite{Golb}, one has

 \begin{equation}\label{Golberg} f_i:=\partial(\det{\cal H}_m)/\partial x_i=\sum_{k+l=i+2} \Delta^k_l
 \end{equation}
where $\Delta ^{k}_l$ is the $(m-1)$-minor of  $\mathcal{H}$ obtained by omitting the $l$th row and $k$th column.
Since $\mathcal{H}_m$ is symmetric,   $\Delta^k_l =\Delta ^l_k$.

On the other hand, from (\ref{GP_trick}), by \cite[Lemma2.3]{GP} and \cite[Corollary2.2]{Conca}, one has

 \begin{equation}\label{delta}
 \Delta ^j_i=\sum_{l=1}^i[1,\cdots,\widehat{l},\cdots,\widehat{(j+i+2-l)},\cdots,n+1]
 \end{equation}
where $j\geq i.$

Collecting the information yields

{\small
\begin{equation}\label{star}
f_j=
\left\lbrace
           \begin{array}{c l}
          \sum_{i=0}^{j/2}(j+1-2i)[1,\cdots,\widehat{(i+1)},\cdots,\widehat{(j+2-i)},\cdots,n+1]     &  j<n,\\
          \\
           \sum_{i=1}^{n-j/2}(2n+1-j-2i)[1,\cdots,\widehat{(i+1+j-n)},\cdots,\widehat{(n+2-i)},\cdots,n+1]   &  j\geq n.
           \end{array}
         \right.
\end{equation}
}

Note that, for any $j$,  the minors appearing in the above expression of $f_j$ are incomparable.
(In the above illustrative diagram minors horizontally aligned are present in one and the same partial derivative.)

\begin{Proposition}\label{integrality}
The radical of the gradient ideal $J$ of a Hankel determinant is the ideal $P$ generated by the $(m-1)$-minors of the matrix.
\end{Proposition}
\demo
By symmetry, it suffices to consider the minors present in one of the partial derivatives  $f_j$, for $j=0,\cdots,2n-2$.
We proceed by descending induction.
For $j=2n-2, 2n-3$ there is nothing to prove since the corresponding partial derivatives are themselves minors, to wit,
$f_{2n-2}=[1,2,\ldots,n-1]$ and $f_{2n-3}=2[1,2,\ldots,n-2,n]$ respectively.

We provide the next inductive step to make the argument clearer.
The next minors in the diagram of partial order are exactly those present in the expression
of $f_{2n-4}$ according to  (\ref{star}).
Consider the Pl\"ucker relation containing the term  $[1,\ldots,n-3,n-1, n][1,\ldots,n-2,n+1]$:
\begin{equation}\label{plucker}
\begin{array}{lcl}
[1,\ldots,n-3,n-1, n][1,\ldots,n-2,n+1] & = & 1/2[1,\ldots,n-3,n-1,n+1]f_{2n-3} \\
 & - & f_{2n-2}[1,\ldots,n-3,n,n+1].
\end{array}
\end{equation}

Now,  $[1,\ldots,n-3,n-1, n]$ satisfies  the following equation:
{\small
\begin{equation}\label{xequation}
X^2 -1/3 X ( 3[1,\ldots,n-3,n-1, n]+[1,\ldots,n-2,n+1] ) + 1/3 [1,\ldots,n-3,n-1, n][1,\ldots,n-2,n+1]=0.
\end{equation}
}

By (\ref{star}) one has $f_{2n-4}=3[1,\ldots,n-3,n-1, n]+[1,\ldots,n-2,n+1]$,  and by (\ref{plucker}) the constant coefficient of
(\ref{xequation}) lies in $J$.
 Hence $[1,\ldots,n-3,n-1, n] \in \sqrt{J}.$

\smallskip

Now suppose that $j<2n-4$ and that, for any $t$ such that $j<t\leq 2n-2$,  any minor present in the expression of $f_t$ belongs to $\sqrt{J}$. Let
 $\Delta:=[1,\ldots,\widehat{i+1+j-n},\ldots,\widehat{n+2-i},\ldots,n+1]$ be a minor present in $f_j$,
 where $1 \leq i \leq [(n-j)/2]$ and $j\geq n$, and write
 $f_j=({\rm coef})\Delta+ \sum$, where $\sum$ denotes the complementary $k$-linear combination of minors as in (\ref{star}).
 Multiplying by $\Delta$ yields $f_j\Delta=({\rm coef})(\Delta)^2+\Delta\,\sum$.

Making everything explicit via (\ref{star}), we see that $\Delta$ satisfies  the following equation:
 \begin{equation}\label{equation}
 X^2 - \frac{1}{2n+1-j-2i} f_j\,X+\frac{1}{2n+1-j-2i} g=0,
 \end{equation}
 where
 {\small
$$g=\sum_{i\neq l=1}^{[(n-j)/2]} [1,\ldots,\widehat{i+1+j-n},\ldots,\widehat{n+2-i},\ldots,n+1] [1,\ldots,\widehat{l+1+j-n},\ldots,\widehat{n+2-l},\ldots,n+1].$$
}

Now, for any $i,l$ such that  $1 \leq i < l \leq [(n-j)/2]$, one has the following  Pl\"ucker relation:
{\scriptsize
\begin{eqnarray}\nonumber
0&=&[1,\ldots,\widehat{i+1+j-n},\ldots,n+2-l,\ldots,\widehat{n+2-i},\ldots,n+1] [1,\ldots,i+1+j-n,\ldots,\widehat{l+1+j-n},\ldots,\widehat{n+2-l},\ldots,n+1]\\\nonumber
&-&{\bf [1,\ldots,{i+1+j-n},\ldots,\widehat{n+2-l},\ldots,\widehat{n+2-i},\ldots,n+1]}
[1,\ldots,\widehat{i+1+j-n},\ldots,\widehat{l+1+j-n},\ldots,{n+2-l},\ldots,n+1]\\ \nonumber
&+&[1,\ldots,\widehat{i+1+j-n},\ldots,\widehat{n+2-l},\ldots,{n+2-i},\ldots,n+1]
{\bf [1,\ldots,i+1+j-n,\ldots,\widehat{l+1+j-n},\ldots,n+2-l,}\\ \nonumber
&&{\bf \ldots,\widehat{n+2-i},\ldots,n+1]}\nonumber
\end{eqnarray}
}

This is nearly a straightening relation as the minors in boldface are less than the starting minor $[1,\ldots,\widehat{i+1+j-n},\ldots,\widehat{n+2-i},\ldots,n+1]$
in the partial order -- but this is all we need. Since these two boldfaced minors are present in
$f_t$ for some $t>j$,  by the inductive hypothesis they belong to $\sqrt{J}$.
 Therefore $g \in \sqrt{J}$. It follows from (\ref{equation}) that $\Delta^2\in \sqrt{J}$, hence $\Delta\in \sqrt{J}$, as required.
 \qed

 \medskip

Drawing upon the preceding results, one has:

\begin{Theorem}\label{p_primary_part}
The minimal component of the primary decomposition of $J$ in $R$ is $P$.
\end{Theorem}
\demo
It suffices to prove that  $P=\sqrt{J}$ and that $P$ is the $P$-primary part of $J$.
The first assertion is the content of Proposition~\ref{integrality}.
As for the second, let $\mathcal{P}$ denote the $P$-primary part of $J$.
Since $\mathcal{P}$ is the contraction of $\mathcal{P}_P$, it suffices to show that $J_P=P_P$, i.e., that $\ell(J_P)=\ell(P_P)$.

By the associativity formula of multiplicities, one has $e(R/J)\geq e(R/P)\ell (R_P/J_P)$, where $\ell$ denotes length.
Since $J'=({\rm in}(f_0),\ldots, {\rm in}(f_{2m-2}))\subset {\rm in}(J)$ also has codimension $3$, using Lemma~\ref{multiplicity_inequality}, yields:
$$e(R/P)\ell(R_P/J_P)\leq e(R/J)=e(R/{\rm in}(J))\leq e(R/J')<2e(R/P).$$
Clearly this is only the case if $\ell(R_P/J_P)=1$.
\qed

\begin{Corollary}\label{colon_by_P}
One has $e(R/J)=e(R/P)\,${\rm ;} in particular, $J:P$ has codimension at least $4$.
\end{Corollary}
\demo The equality of multiplicities is clear from Proposition~\ref{p_primary_part} by the associativity formula.
From the exact sequence $0\rar P/J\rar R/J\rar R/P\rar 0$ we then conclude that
$$\dim P/J<\dim R/J=2m-1-3.$$
Since $J:P$ is the annihilator of $P/J$ we are through with the subsumed statement.
\qed

\bigskip

We next state three open questions.
Due to large computer verification, natural interspersing of the various properties so far and last, but not least, aesthetic formulation, we risk to call them conjectures.
They are listed in increasing order of difficulty (to our best guessing).

\begin{Conjecture}\label{reduction}\rm
Let $\mathcal{H}_m$ denote the Hankel matrix of size $m\times m$.
Let as above $J$ denote the gradient ideal of its determinant and $P=I_{m-1}(H)$.
Then, for $0\leq i\leq m-2$ one has
$$JP^i:P^{i+1}=I_{m-2-i}(H).$$
\end{Conjecture}

\begin{Corollary} {\rm (of Conjecture~\ref{reduction})}
\begin{enumerate}
\item[{\rm (i)}] $Q:=I_{m-2}(H)$ is the only other associated prime of $R/J$ {\rm (}necessarily embedded{\rm )}
\item[{\rm (ii)}] $JP^{m-3}:P^{m-2}=(x_0,\ldots,x_{2m-2})$
\item[{\rm (iii)}] $JP^{m-2}:P^{m-1}=(1)$, i.e., ${\rm red}_J(P)=m-2\,${\rm ;} in particular, $J$ is a minimal reduction of $P$.
\item[{\rm (iv)}] The partial derivatives of $\det(H)$ are algebraically independent over $k$.
\end{enumerate}
\end{Corollary}
\demo (Highlights)

(i) Use Theorem~\ref{p_primary_part}. Note that this item  is now equivalent to showing that $J:P=Q$.
Indeed, let $J=P\cap {\cal N}$, stand for the primary decomposition of $J$, where ${\cal N}$ denotes the intersection of the remaining (necessarily embedded) primary components.
Then $Q=J:P={\cal N}:P$, hence ${\cal N}$ is $Q$-primary.
It follows that $Q$ is the only associated embedded prime of $R/J$.

(ii)  This guarantees that the reduction number is not smaller than $m-2$.
In this regard we risk to say, moreover, that $\mu(P^{m-2})=\mu(JP^{m-3})+1$, so (ii) would correspond to having the entries in the (say) last coordinate of the syzygies of $P^{m-2}$ generate the maximal ideal.
Thus, previous information about the syzygies of the powers of $P$ may be useful.

(iii) Clearly, $J$ is a reduction of $P$. To see that it is a minimal reduction it suffices to see that it is minimally generated by $2m-1$ forms, while this is the analytic spread of $P$ because the dimension of the $k$-subalgebra of $R$ generated by the $(m-1)$-minors is $2m-1$.

(iv) This follows since $J$ is a minimal reduction, hence is generated by analytically independent forms.
But, since for forms of the same degree, ``analytical independence'' and ``algebraic independence'' are the same notion,
we get that the partials themselves are algebraically independent.
Equivalently, the inclusion of the algebras $k[J_{n-1}]\subset k[P_{n-1}]$ is a Noether normalization.
\qed

\begin{Conjecture}\label{linear_type}\rm
 $J$ is an ideal of of linear type.
\end{Conjecture}

\begin{Conjecture}\label{non-homaloidal}\rm
The determinant of ${\cal H}_m$ is not homaloidal.
\end{Conjecture}

\begin{Remark}\rm
(i) It is shown in \cite{M_thesis} that the syzygy matrix of $J$ has linear rank $3$.
Together with an affirmative answer to Conjecture~\ref{linear_type} this would imply that $\det {\cal H}$ is {\em not} homaloidal (last of the three conjectures above).
In the subsequent part we will prove both conjectures in the case $m=3$.

(ii) Regarding assertion (iv) of the above conjectured corollary, it is indeed true that the Hankel partial derivatives are algebraically independent over $k$ - a proof is given in \cite[Proposition 3.3.11]{M_thesis} in terms of the non-vanishing of the corresponding Hessian determinant. In this regard it would be pertinent to give a proof of the parabolism of the Hankel determinant $f$, i.e., that $f$ divides the Hessian $H(f)$ and to find the exact multiplicity exponent.
Note that this exponent is $1$ in the case of $m=3$ (Example~\ref{parabolic_hankel}).
\end{Remark}

\subsection{Cases study}

This part is devoted to examining small sizes of catalecticants.
The reason for this section is either because for small value of $m$ one has special properties inexistent for larger size, or else because some behavior is conjectured to hold in any size and yet we lack the tools to write a complete proof.

\subsubsection{Hankel $3\times 3$}

\smallskip

This is the matrix
$$\mathcal{H}_3=\left(
  \begin{array}{ccc}
    x_0 & x_1 & x_2 \\
    x_1 & x_2 & x_3 \\
    x_2 & x_3 & x_4 \\
  \end{array}
\right)
$$

\begin{Proposition}\label{hankel3x3}
Let ${\cal H}_3$ denote the generic $3\times 3$ Hankel matrix in the
variables $x_0,\ldots,x_4$ and let $J\subset R:=k[x_0,\ldots,x_4]$ denote the
gradient ideal of $\det {\cal H}_3$.
\begin{enumerate}
\item[{\rm (a)}] Let $P:=I_2({\cal H}_3)\subset R\,${\rm ;} then $P$ and $\fm:=(x_0,\ldots,x_4)$
are the only associated primes of $R/J$.
\item[{\rm (b)}]  $\det {\cal H}_3$ is not homaloidal.
\end{enumerate}
\end{Proposition}
\demo
Drawing on the previous results (Theorem~\ref{p_primary_part}), we only have to show that $\fm$ is an associated prime.
For this it is enough to show that $\fm\subset J:P$.
If the lower bound for $J:P$ from Corollary~\ref{colon_by_P} could be improved by $1$, we would be done since at any rate $J:P$ is contained in some associated prime of $R/J$.

Instead we argue directly within the details of this case.
We switch freely back and forth between the Hankel matrix and the associated $2\times 4$ Hankel matrix (\ref{GP_trick}), where the $2$-minors can be thought as maximal minors.
By the shape of the generators of $J$, it suffices to show that, for $i=0,\ldots,4$, either $x_i\Delta_{23}\in J$ or $x_i\Delta_{14}\in J$.
Taking the $3\times 2$ submatrix of the given $3\times 3$ Hankel matrix consisting of the last two columns,
we immediately get that $x_3\Delta_{23}, x_4\Delta_{23}\in J$. By a similar token, the $2\times 3$ submatrix
formed by the first two rows yield $x_0\Delta_{23}, x_1\Delta_{23}\in J$.
To deal with $x_2$ choose $\Delta_{14}$ instead. This time around, consider the $2\times 3$ submatrix of the associated $2\times 4$ Hankel matrix (\ref{GP_trick})
 consisting of the first and last two columns.
Then $x_2\Delta_{14}\in J$, as required.

\medskip

(b) By (a), the saturated ideal $J:\fm^{\infty}$ coincides with $P$.
Supposing that $J$ defines a Cremona map, one would have that the initial degree of $P/J$ is at least $2+1=3$ (\cite[Proposition 1.2]{PanRusso}).
But this is nonsense since $P$ admits a generator of degree $2$ (minor) that does not belong to $J$.
\qed

\medskip

Note that part (b) above proves Conjecture~\ref{non-homaloidal} for $m=3$, while the proof of part (a) showed that $J:P=\fm=I_1({\cal H}_3)$, which is half the statement of Conjecture~\ref{reduction}.
To complete the other half, we need to show that $JP:P^2=(1)$.
For convenience we isolate this result as a proposition.

\begin{Proposition}
Let ${\cal H}_3$ denote the generic $3\times 3$ Hankel matrix, let $J\subset R$ denote the
gradient ideal of $\det {\cal H}_3$ and $P:=I_2({\cal H}_3)$.
Then $P^2\subset JP\,${\rm ;} in particular, $J$ is a minimal reduction of $P$ with reduction number $1$.
\end{Proposition}
\demo
Since $P=(J,\Delta)$, where $\Delta=\Delta_{23}=x_2^2-x_1x_3$, it suffices to show that $\Delta^2\in JP$.
As in the proof of Theorem~\ref{p_primary_part} a direct calculation up to nonzero coefficients and adjusting signs yields the relation
$$\Delta^2=f_2\Delta + \Delta_{14}\Delta_{23}.$$
Clearly, $f_2\Delta\in JP$.
Now, using the Pl\"ucker relation
$\Delta_{12}\Delta_{34}-\Delta_{13}\Delta_{24}+\Delta_{14}\Delta_{23}=0$,we get that $\Delta_{14}\Delta_{23}\in JP$ as well.
\qed

\medskip

The conjecture about the linear type property is proved in the following theorem.

\begin{Theorem}
Let ${\cal H}_3$ denote the generic $3\times 3$ Hankel matrix and let $J\subset R$ denote the
gradient ideal of $\det {\cal H}_3$.
Then $J$ is an ideal of of linear type.
\end{Theorem}
\demo To show the linear type property we use the criterion of \cite[Theorem 9.1]{Trento}, informally
known as the ``$(F_1)$ $+$ sliding depth'' criterion.
Now, $(F_1)$ (also known as $(G_{\infty}$) is the following:

\smallskip

{\sc Claim}: $\mu(J_Q)\leq \hht Q$, for every prime ideal $Q$.

\smallskip

To prove this assertion we have to drill quite a bit through the syzygies of the generators of $J$.
Since $P$ has the expected codimension ($=3$), its presentation matrix is linear, as a piece of the Eagon--Northcott
complex for the maximal minors of the Hankel matrix

\begin{equation}\label{GP3}
\left( \begin{matrix}
x_0&x_1&x_2&x_3\\
x_1&x_2 & x_3 &x_4\\
\end{matrix} \right)
\end{equation}

Noting that the syzygies derive from the Laplace relations coming from the four $2\times 3$ submatrices
of (\ref{GP3}), and ordering the minors to adjust to the order in which they appear along the usual ordering of the partial
derivatives, one obtains the following matrix:

$$ \mathcal{S}=\left( \begin{matrix}
0&0&0&x_0&-x_1&0&x_1&x_2\\
0&x_0 &0&0&x_2&x_1&0&-x_3\\
0&-x_1&0&-x_2&0&-x_2&-x_3&0\\
x_0&0&x_1&0&-x_3&0&0&x_4\\
-x_1&0&-x_2&x_3&0&0&x_4&0\\
x_2 & x_3 & x_3 &0 &0 & x_4 & 0&0\\
\end{matrix} \right )$$
Let $\mathcal{S}_j\, (1\leq j\leq 8)$ denote the $j$th column of this matrix.
A linear syzygy of the ordered partial derivatives comes by ``tweaking'' a $k$-linear combination of columns of $\mathcal{S}$ satisfying
the following restrictions: first, the entry on the fourth coordinate is $3$ times the one on the third coordinate; second,
the set of entries along the $3$rd row of the involved columns is the same as the set of entries along the $4$th row.
After a detailed inspection of these restrictions, we get at least the following possibilities of pairs: $\{\mathcal{S}_2,\mathcal{S}_3\}$, $\{\mathcal{S}_4,\mathcal{S}_6\}$ and
$\{\mathcal{S}_5,\mathcal{S}_7\}$ and the resulting syzygies being
$3\mathcal{S}_3-\mathcal{S}_2, -3\mathcal{S}_5-\mathcal{S}_7,
-\mathcal{S}_6+\mathcal{S}_4$.

In this way the ordered partial derivatives admit the following linear syzygies
\begin{equation}\label{linear_syzygies3}
\left( \begin{matrix}
0   &  -2x_0 & 4x_1\\
x_0 & -x_1 & 3x_2\\
2x_1 &  0  & 2x_3\\
3x_2 & x_3 & x_4\\
4x_3 & 2x_4  &  0
\end{matrix} \right ).
\end{equation}
yielding a matrix of rank $3$.

We now contend that, together with the Koszul syzygies, they suffice to check the stated property of $J$.
For this, we note that this property is equivalent to a property of the Fitting ideals of $J$, and can be
stated as follows:
$$\hht I_t(\phi)\geq \rk (\phi)-t +2=4-t+2=6-t, \;\; {\rm for}\;\; 1\leq t\leq 4.$$
Obviously, it suffices to show these estimates for the above submatrix $\psi$.
For $t=1,2$ the linear syzygies so far described immediately give the required inequalities.
For $t=3,4$, we use the matrix of Koszul syzygies where one knows that the codimension of any of its Fitting ideals is
$\geq \hht (J)=3$.
So much for the property $(F_1)$.

\smallskip

The sliding depth condition reads as

\smallskip

{\sc Claim}: $\depth(H_i)_Q \geq \hht(Q) - \mu(J_Q) + i$, for every prime  $Q \supset J$
with $0 \leq i \leq\mu(J_Q) -\hht(J_Q )$.

Here $H_i$ denotes the $i$th Koszul homology module on the generators of $J$

\smallskip

Note that if $Q\supset J$ then necessarily $Q\supset P$.
We first assume that $Q=\fm$, in which case $\mu(J_{\fm})=\mu(J)=5$, $\hht(J_{\fm})=\hht(J)=3$, hence
$\mu(J_{\fm})-\hht(J_{\fm})\leq 2$.
In this case, our range for $i$ is $0\leq i\leq 2$.
For $i=2$ the inequality is automatic since $H_2$ is Cohen--Macaulay and it is trivial if $i=0$.
For $i=1$  we need to show that $H_1$ has depth at least $1$.
This computation has been carried out with {\em Macaulay} by showing that the homological dimension of $H_1$ over $R$ is $4$.

Next assume that $Q\subsetneq \fm$, hence $\mu(J_Q)\leq\hht(Q)\leq 4$ (by the first claim above).
Clearly, $\hht(J_Q)=\hht(P_Q)=3$. Since $P$ is the minimal part of $J$ and $\fm$ is the only other associated prime of $R/J$, one has $J_Q=P_Q$. Therefore, $\mu(J_Q)-\hht(J_Q)=\mu(P_Q)-\hht(P_Q)=3-3=0$
since $R/P$ is an isolated singularity.
Then $i=0$ is the only case. One has $\depth (H_0)_Q=\depth (R_Q/P_Q)=2 = \hht(Q)-\mu(J_Q) +0$.

\smallskip

Thus, $J$ is an ideal of linear type.
\qed

\begin{Remark}\rm
The linear syzygies (\ref{linear_syzygies3}) in fact generate all linear syzygies of $J$ as is shown in \cite[Theorem 3.3.5]{M_thesis} in greater generality. Together with the linear type property, it implies via Proposition~\ref{polar_and_lineartype} that the Hankel determinant of size $3$ is not homaloidal.
\end{Remark}

\subsubsection{Catalecticant: $2$-leap $3\times 3$}

In contrast with the Hankel case, we don't have as yet a general theory of the associated primes of the gradient ideal of the catalecticant determinant.

In the case of a $3\times 3$ catalecticant, since the upper extreme case of the leap ($r=3$) is the generic $3\times 3$ matrix and the lower extreme is the Hankel matrix, both of which have been taken care of, we are left with the case  $r=m-1=2$.

This is the matrix
$$\mathcal{C}=\mathcal{C}_{3,2}=\left(
  \begin{array}{ccc}
    x_0 & x_1 & x_2 \\
    x_2 & x_3 & x_4 \\
    x_4 & x_5 & x_6 \\
  \end{array}
\right)
$$

We call attention to the sharp difference at the outset between this  catalecticant case and the previous $3\times 3$ Hankel matrix:
here the ideal of $2$-minors of $\mathcal{C}$ is radical, but not prime. Overall the structure of the primary components of $R/J$ is a lot more intricate as we now explain.

\begin{Proposition}\label{2-cat3x3}
Let $J\subset R:=k[x_0,\ldots,x_6]$ denote the
gradient ideal of $\det \mathcal{C}$.
Consider the associated  $2$-leap catalecticant
$$\widetilde{\mathcal{C}}:=
\left(
  \begin{array}{ccccc}
    x_0 & x_1 & x_2 & x_3 & x_4 \\
    x_2 & x_3 & x_4 & x_5 & x_6 \\
  \end{array}
\right)
$$
\begin{enumerate}
\item[{\rm (a)}] Let $I:=I_2(\mathcal{C})\subset R\,$ and $P:=I_2(\widetilde{\mathcal{C}})${\rm ;}
then the minimal part of $J$ is the radical ideal $I=P\cap (I:P)=P\cap (x_0,x_2,x_4,x_6)$ and
$Q:=(I:P)+P=(x_0,x_2,x_4,x_6, x_1x_5-x_3^2)$
is the only embedded prime of $R/J$.
\item[{\rm (b)}]  The generators of $J$ are algebraically independent over $k$ and the syzygy matrix
of $J$ has maximal linear rank. In particular, $\det \mathcal{C}$ is homaloidal.
\end{enumerate}
\end{Proposition}
\demo
Since the data are so very explicit and of small size, one can set up a computation in {\sc Macaulay\/} \cite{Macaulay} to
check all the assertions.
We choose instead to play the theory against the background as much as possible before resorting to a computer calculation.

(a) First, the present situation confronts two matrices by adapting a known general principle (\cite[Lemme 2.3]{GP}), which in this case tells that the $2$-minors
of $\mathcal{C}$ are the $2$-minors of $\widetilde{\mathcal{C}}$ excluding the minor $x_1x_5-x_3^2$
of columns $2$ and $4$ of the latter.
(Of course this can also be checked directly by inspection.)

Thus, in terms of the stated notation, one has $P=(I,b)$, where $b=:x_1x_5-x_3^2$.
Since $P$ is a prime ideal of codimension $4$, $I$ has codimension at least $3$.
But $I\subset (x_0,x_2,x_4,x_6)$, the latter being a prime ideal of codimension $4$ and, clearly, a minimal prime
thereof.
Therefore, these two prime ideals are minimal primes of $R/I$.
We claim that their intersection is contained in $I$.
To see this, it suffices to show that $b\cdot(x_0,x_2,x_4,x_6)\subset I$ as $b$ is a non-zero-divisor modulo
$(x_0,x_2,x_4,x_6)$.
But this is the content of Laplace rule as applied to suitable $2\times 3$ submatrices of
$\widetilde{\mathcal{C}}$ involving columns $2$ and $4$.
(One can similarly show that $I:P=(x_0,x_2,x_4,x_6)$ as written on the statement, but we will have no use for it.)

We now deal with the gradient ideal $J$.
In particular, we can express the partial derivatives of $f:=\det \mathcal{C}$ in terms of the $2$-minors of
$\widetilde{\mathcal{C}}$ rather than those of $\mathcal{C}$ itself, whereby $x_1x_5-x_3^2$ will
not show.
Letting $\Delta_{ij}$ denote the $2$-minor of $\widetilde{\mathcal{C}}$ with $i$th and $j$th columns ($i<j$), one has
{\small
\begin{equation}\label{gens_of_J_cat}
\{f_0=\Delta_{45}, f_1=-\Delta_{35}, f_2=  2\Delta_{34}-\Delta_{25}, f_3=\Delta_{15}, f_4=2\Delta_{23}-\Delta_{14}, f_5=-\Delta_{13},
f_6=\Delta_{12}\}
\end{equation}
}
We first check that $J$ has codimension $4$.
Namely, the sequence $f_0,f_2, f_5,f_6$ is a regular sequence.
In order to see this the easiest is perhaps to note that the three minors $f_0, f_5,f_6$ are known to form a regular sequence by specializing from the generic case (with $(x_2,x_3,x_4)$
being a minimal prime of minimal codimension).
Then show that $f_2$ is a nonzerodivisor thereof --  to get around this point we resort to the computer since the structure of the associated primes of $R/(f_0,f_5,f_6)$ is rather involved.

As a result, $P$ and $I:P$ are minimal primes of maximal dimension of $R/J$, hence $I=P\cap (I:P)$ contains the unmixed part $J^{\rm un}$ of $J$.
 In order to show that  $J^{\rm un}=I$ one would first show that  $R/J$ has no minimal primes other than $P$ and $I:P$, then that  $J^{\rm un}$ is a radical ideal.
 The multiplicity of $R/J$ is at least $e(R/P)+e(R/(I:P))=5+1$,
 where the second of these is obvious and the first follows from the fact that the Eagon--Northcott resolution of $R/P$ is pure and $2$-linear.
 However, showing before hand that $e(R/J)\leq 6$ seems like a tough matter.
 We found no easy hand calculation for this part, not to mention the  resources one had in the Hankel case that we lack here.
 As a dissonant note, here $I^2\subset J$ but $I^2\not\subset JI$ which dashes our hope for any reduction theory to come in.
 Resorting to {\em Macaulay} one easily computes the equality $J^{\rm un}=I$. As a result, $e(R/J)=6$ and hence, using the exact sequence $0\rar I/J\lar R/J \lar R/I \rar 0$, where $R/J$ and $R/I$ have same dimension and multiplicity, one derives that the annihilator $R/J:I$ has codimension at least $5$. This completes the proof that $R/J$ has no other minimal primes and any embedded prime has codimension at least $5$.
 At this point, an additional computation shows that $J:I=(x_0,x_2,x_4,x_6,x_1x_5-x_3^2)$ and this is then the only embedded prime of $R/J$.

(b) When char$(k)=0$, the partial derivatives are algebraically independent over $k$ if and only if
the Hessian determinant does not vanish (both express the condition that the polar map has image
of maximal dimension).
If one is in peace by accepting the result in Example~\ref{parabolic_exs}, the Hessian determinant is clearly nonzero. Since the non-vanishing is much less precise than that result, we might choose to proceed directly, as follows.
First, a tedious but straightforward calculation gives the Hessian matrix:

$$\left( \begin{matrix}
0  &   0  &   0  &   x_6 &  -x_5 & -x_4 & x_3\\
0  &   0  &   -x_6 & 0  &   2x_4 & 0  &   -x_2\\
0   & -x_6 & 2x_5 & -x_4 & -x_3 & 2x_2 & -x_1\\
x_6 &  0  &   -x_4 & 0  &   -x_2 & 0  &   x_0\\
-x_5 & 2x_4 & -x_3 & -x_2 & 2x_1 & -x_0 & 0\\
-x_4 & 0  &   2x_2 & 0  &   -x_0 & 0  &   0\\
x_3 &  -x_2 & -x_1 & x_0 &  0 &    0  &   0
\end{matrix} \right)
$$
Evaluating at $(0,0,1,0,0, 1, 1)\in \mathbb{Q}^7\subset k^7$ yields the numerical matrix
$$
\left( \begin{matrix}
0 & 0 & 0 & 1 & -1 &0 & 0\\
0 & 0 & -1 & 0 & 0 & 0 &-1\\
0 & -1& 2 & 0 & 0 & 2 & 0\\
1 & 0 & 0 & 0 & -1 & 0 & 0\\
-1 & 0 & 0 & -1 & 0 & 0 & 0\\
0 &  0 & 2 & 0 & 0 & 0 &0\\
0 & -1& 0 & 0 & 0 & 0 & 0
\end{matrix} \right)
$$
whose determinant is immediately computed by Laplace: its value is $8\in \mathbb{Q}$.
Therefore the numerical matrix is non-singular and hence the Hessian determinant does not vanish.

We next show that the linear rank of the partials is maximal, i.e.,  $6$.
Together they imply by Proposition~\ref{polar_and_lineartype} that $f$ is homaloidal.

Clearly, from the form of the partial derivatives in (\ref{gens_of_J_cat}), any
linear syzygy of those is a  linear syzygy of the ideal $I$ of $2$-minors of $\mathcal{C}$.
More precisely, order the $2$-minors of so as to adjust to the order in which they
appear along the partials in (\ref{gens_of_J_cat}). Then a linear syzygy $(l_0,\ldots,l_8)^t$ of
this ordered minors gives rise to one of (\ref{gens_of_J_cat})
if and only $l_2=-2l_3$ and $l_5=-2l_6$ in which case the former will have coordinates
$(l_0,-l_1,-2l_3,l_4,-2l_6,-l_7, -l_8)^t$.

In order to search for such syzygies of $I$ one resorts to its presentation.
We note that though $I$ is actually a nice Gorenstein ideal -- being a specialization of the generic case
-- it is easier to work with $P$ since its presentation is a piece of the Eagon--Northcott
complex, whereby the syzygies derive from the Laplace relations coming from the five $2\times 3$ submatrices
of $\widetilde{\mathcal{C}}$. Thus, with the previous order of the generators of $I$, leaving out
the syzygies of $P$ that involve the extra minor $x_1x_5-x_3^2$ with a nonzero coefficient, one obtains the following matrix:
{\small
$$
\left( \begin{array}{cccccccccccccccc}
0&     0&     0&     0&     x_0 &  0&     x_4 &  0&     0&     0&     0&     x_2 &  0  & -x_2 & 0&  x_1\\
0&     0&     0&     x_0 &  0&     0&     -x_5 & x_1 &  0&     0&     x_2 &  0&  -x_3 & x_3 &  0&  0\\
0&     0&     x_0 &  0&     0&     0&     x_6 &  0&     x_2 &  0&     0&     0&  0 & -x_4 & x_1 &  x_3\\
0&     x_0 &  0&     0&     0&     0&     0&     -x_2 & 0&     x_2 &  0&     0&  x_4  & 0&  0& -x_3 \\
0&     -x_1 & 0&     -x_2 & -x_3 & 0&     0&     0&     0&     -x_3 & -x_4 & -x_5 & 0 & 0& 0& 0\\
x_0 &  0&     0&     0&     0&     x_2 &  0&     x_4 &  0&     0&     0&     0&  -x_6 & 0&   x_3 &  x_5 \\
0&     0&     -x_2 & 0&     x_4 &  0&     0&     0&     -x_4 & 0&     0&     x_6 &  0& 0&   -x_3 & 0\\
-x_1 & 0&     x_3 &  x_4 &  0&     -x_3 & 0&     0&     x_5 &  0&     x_6 &  0&  0 & 0&  0&  0\\
x_2 &  x_4 &  0&     0&     0&     x_4 &  0&     0&     0&     x_6 &  0&     0&  0 &0&  x_5 &  0
\end{array} \right)
$$
}
where the last two syzygies cannot be expressed as Laplace relations and compensate for the
omission of the Laplace relations coming from the $2\times 3$ submatrices of $\widetilde{\mathcal{C}}$
fixing columns $2$ and $4$.

Let $\mathcal{S}_j\, (1\leq j\leq 16)$ denote the $j$th column of this matrix.
An additional restriction is that if a $k$-linear combination of a set of these columns
 is to produce a syzygy $(l_0,\ldots,l_8)^t$ such that $l_2=-2l_3$ and $l_5=-2l_6$
 then the set of entries along the $3$rd row of the involved columns in the combination has to match the set of entries
along the $4$th row, and the same regarding the $6$th and $7$th rows.
By inspection, we get at least the following sets
$$\{\mathcal{S}_2,\mathcal{S}_3,\mathcal{S}_6\}, \{\mathcal{S}_8,\mathcal{S}_9,\mathcal{S}_{10}\},
\{\mathcal{S}_{12},\mathcal{S}_{13},\mathcal{S}_{14}\}, \{\mathcal{S}_5,\mathcal{S}_8,\mathcal{S}_9,\mathcal{S}_{10}\}.
$$
Each one of these sets gives rise to a linear syzygy of $I$ satisfying the required condition, from
which we derive as explained the corresponding linear syzygy of $J$.
This gives $4$ linear syzygies of $J$.
Furthermore, we have two additional linear syzygies coming from the Laplace rule as applied
to the generators $f_1,f_3,f_5$ formed with the corresponding $2\times 3$ submatrix of
$\widetilde{\mathcal{C}}$.
Therefore, the following linear syzygies are obtained:
$$
\left( \begin{array}{cccccc}
0 & 0&       0 &         0 &         x_2 &      x_0 \\
x_0& x_2 &   0 &       x_1 &       2/3x_3 &    2/3x_1 \\
0 & 0&       x_0 &  -1/4x_2 &     2/3x_4 &     1/3x_2 \\
x_2 & x_4 &  x_1 &   3/4x_3 &    1/3x_5 &         0\\
0 & 0&       2x_2&  -1/2x_4 &    1/3x_6 &     -1/3x_4 \\
x_4 & x_6&   2x_3&   1/2x_5 &         0 &     -2/3x_5 \\
0 & 0&       3x_4&  -3/4x_6 &         0 &      -x_6
\end{array} \right)
$$
Evaluating the upper $6\times 6$ submatrix  at $(1,0,0,1,1,0,0)$ gives the following matrix
$$
\left( \begin{array}{cccccc}
0 & 0 & 0 & 1& 0 & 0\\
0 & 1 0 & 0 & 0 & 2/3\\
0 & 0 & 1 0 & 0 & 2/3\\
1 0 & 0 & 0 & 3/4 &  0\\
0 & 0 & 0 & -1/3 & -1/2 & 0\\
0 & 1 & 2 & 0 & 0 & 0
\end{array} \right)
$$
whose determinant is easily calculated by Laplace and has value $-1$.
Therefore, the rank of the above matrix of linear syzygies is $6$.
\qed

\begin{Remark}\label{dual_variety}\rm
 A computation with {\em Macaulay} (\cite{Macaulay}) actually shows that $J$ is of linear type and the inverse map has degree $6$
(maximum possible in this situation). Therefore, the inversion factor has degree $11$ and the Hessian determinant
is a factor thereof, with complementary factor the equation $g$ of the dual variety to the twisted cubic as found in Example~\ref{parabolic_exs}.
This fact clamors for a geometric explanation.
\end{Remark}

\begin{Question}\rm For what values of $1\leq r\leq m$ is the determinant $f_{m,r}$ of the generic $r$-leap catalecticant
 homaloidal?
\end{Question}

We believe that the answer is: $f_{m,r}$ is homaloidal if and only if $r=m$ or $r=m-1$.
Next is some computational and theoretical evidence.

\medskip

\subsubsection {Additional cases}

Consider the $3$-leap $4\times 4$ catalecticant

\smallskip

$$\mathcal{C}=\mathcal{C}_{4,3}=\left(
  \begin{array}{cccc}
    x_0 & x_1 & x_2 & x_3\\
    x_3 & x_4 & x_5  & x_6\\
    x_6 & x_7 & x_8  & x_9\\
    x_9 & x_{10} & x_{11} & x_{12}
  \end{array}.
\right)
$$

We sketch the main features of this case, including the fact that $\det \mathcal{C}$ is homaloidal.

Let $J\subset R:=k[x_0,\ldots,x_{12}]$ denote the
gradient ideal of $\det \mathcal{C}$.
Then, by
an argument that mixes computational and conceptual results,  the following hold.

\begin{enumerate}
\item[{\rm (a)}] (Ideal structure) Let $I:=I_3(\mathcal{C})\subset R\,$ and $P:=I_3(\widetilde{\mathcal{C}})$, where
$$\widetilde{\mathcal{C}}:=
\left(
  \begin{array}{ccccccc}
    x_0 & x_1 & x_2 & x_3 & x_4 & x_5 & x_6\\
    x_3 & x_4 & x_5 & x_6 & x_7 & x_8 & x_9\\
    x_6 & x_7 & x_8 & x_9 & x_{10} & x_{11} & x_{12}
  \end{array}
\right).
$$
Then the unmixed part of $J$ is the codimension $4$ prime ideal $I$ and moreover, $P$ is a codimension $5$ embedded prime
 of $R/J$ with $I=J:P$.

 The ideal $I$ has codimension $4$ since it specializes from the generic case -- and is, in fact, Gorenstein.
 The matrix $\widetilde{\mathcal{C}}$ is obtained by analogy to the Gruson--Peskine procedure, but of course $I$ and $P$
 differ extensively in both in number of generators as in codimension; both ideals are however prime (cf. the case $m=3$,
 where $I$ was not prime) and are the associated primes of $R/J$.
 A computation with {\em Macaulay} (\cite{Macaulay}) gives that the unmixed part of $J$ is $I$ and the latter is $J:P$.

\item[{\rm (b)}] (Structure of the partial derivatives)
By inspection, $10$ out of the $13$ partial derivatives are $3$-minors  of $\mathcal{C}$ up to signs.

Of these, $8$ are up to signs the  $3$-minors of the two submatrices
of $\mathcal{C}$ consisting of columns $1,2,4$ and $1,3,4$, respectively.
Each four generate a prime ideal of height $2$, hence the sum has height $\geq 3$.
But then it has actually  codimension $3$ since it is clearly contained in  $I_2(\psi)$, where $\psi$ denotes the Hankel matrix formed with the first and the last columns of $\mathcal{C}$.
Denote by $J'\subset J$ generated by these $8$ minors.
A guess, inspired on the case of the general $4\times 4$ matrix, was that its unmixed part is $I_2(\psi)$.
A computation with {\em Macaulay} confirmed this and, in addition, gave the equality  $J':I_2(\psi)=I$, which was also expected to hold by specialization from the $4\times 4$ general matrix.

\item[{\rm (c)}] (Linear syzygies)
The syzygy matrix of $J$ has linear rank $11$.

The above ideal $J'\subset J$ has $6$ independent linear syzygies by applying the Laplace trick.
Extending these syzygies to syzygies of all partial derivatives by filling the remaining slots with zeros yields a submatrix of rank $6$ of the syzygy matrix of $J$.
Five additional independent linear syzygies are found by the method of Proposition~\ref{2-cat3x3} and Proposition~\ref{hankel3x3}, bringing up a linear submatrix of rank $11$.
Of course, this can also be computed with {\em Macaulay}.

\item[{\rm (d)}] (Syzygetic part) The syzygetic part $\delta(J)$ has $4$ minimal generators (in standard degree $4$).

The syzygetic part is the kernel $\delta(J)$ of the natural map $H_1\rar R^{13}/JR^{13}$ induced by applying $\otimes_RR/J$
to the presentation sequence $0\rar Z\lar R^{13}\lar J\rar 0$ of $J$; here $H_1$ denotes the first Koszul homology module on the partial derivatives. More explicitly, $\delta(J)=Z\cap JR^{13}/B$, where $B$
denotes the image of the Koszul map on those generators.
Note that $H_1$ is a graded $R$-module with the standard grading of $R$ and that $\delta(J)$ inherits this grading.
Write $\delta(J)=\oplus_{t\geq 1} \delta(J)_{t+3}$
By \cite[Section 1]{SV} one knows that the minimal number of generators of the $R$-module $\delta(J)_{t+3}$ coincides with the
cardinality of a $k$-basis of $\mathcal{J}_{(t,2)}$ consisting of elements of bidegree $(t,2)$, where
$\mathcal{J}$ denotes the bigraded presentation ideal of the Rees algebra of $J$ over $R$.

In particular, the minimal generators of $\delta(J)$, being all of standard degree $4$, give a basis $Q$ of the $k$-vector space of Rees equations of bidegree $(1,2)$.

\item[{\rm (e)}] (Conclusion)
The partial derivatives are algebraically independent over $k$ and $\det \mathcal{C}$ is homaloidal.

Consider the set consisting of the forms of $\mathcal{J}$ of bidegree $(1,1)$ (linear syzygies of $J$) and the set $Q$ above of forms  of bidegree $(1,2)$.

One would claim that the rank of the Jacobian matrix of this set with respect to the variables $\mathbf{x}$ is $12$.
The idea is to look at this matrix in a suitable way to build up a non-zero $12\times 12$ minor thereof involving at least one of the
last $4$ rows with quadratic entries in $\yy$.
This is an additional technical challenge; alternatively,
a computation with {\em Macaulay} gives the expected rank.
By the general criterion of \cite[Theorem 2.18]{AHA}, $\det \mathcal{C}$ is homaloidal.

\end{enumerate}

\medskip

For the  $2$-leap $4\times 4$ catalecticant

$$\mathcal{C}_{4,2}=\left(
  \begin{array}{cccc}
    x_0 & x_1 & x_2 & x_3\\
    x_2 & x_3 & x_4  & x_5\\
    x_4 & x_5 & x_6  & x_7\\
    x_6 & x_{7} & x_{8} & x_{9}
  \end{array}
\right)
$$
we suspect that this case is not homaloidal.
But alas, the sole evidence is computational:
the linear rank is $6$, hence $3$ less than the maximum possible ($=9$). In addition, the Rees algebra has two generators of bidegree $(1,2)$ which come from computing the syzygetic part of the ideal. In principle, it could have additional equations of bidegree $(1,s)$ with $s>2$, although a direct (stalled) computation seems to give no new generators, so the Jacobian dual matrix would not have maximal rank.

\begin{Question}\rm
Is there an explicit formula for the linear rank of the gradient ideal in the case of catalecticants which depends solely on the size and leap?
\end{Question}

\section{Matrices. II: degenerations}

As explained in the Introduction, this section is devoted to understanding the
behavior of certain degenerations of the examples considered in the precious section.
We chose to start with the Hankel degeneration since its homaloidal behavior is fully understood
in \cite[Section 4.1]{CRS}, so we have nothing to add in this respect.
Not so however its homological contents, so to say. In fact we will solve affirmatively
the conjectural note stated in \cite[Remark 4.6 (c)]{CRS}, a step that, as mentioned
there has strong implication to some of the technology developed.
The actual remark there points to a reference that in fact never appeared and, somehow, can be
thought of as being ``replaced'' by the present work.

\subsection{Degeneration of Hankel matrices and their homology}\label{deg_of_cat}

In this part we deal with the determinant of a so-called
sub-Hankel matrix,  considered in \cite[Section 4]{CRS}.
Since a germ of these considerations existed prior to the latter paper, now used amongst the preliminaries
of the present work, we will repeat some of its
relevant aspects here. For consistency and appropriate referencing,
we will keep  the same notation as in \cite{CRS}.

Let $x_0,\ldots,x_n$ be variables over a field $k$ and set

$$
M^{(n)}=M^{(n)}(x_0,\ldots,x_n)= \left(
\begin{matrix}
x_0&x_1&x_2&\ldots &x_{n-2}&x_{n-1}\\
x_1&x_2&x_3&\ldots &x_{n-1}&x_n\\
x_2&x_3&x_4&\ldots &x_n&0\\
\vdots &\vdots &\vdots &\ldots&\vdots &\vdots \\
x_{n-2}&x_{n-1}&x_n&\ldots &0&0\\
x_{n-1}&x_n&0&\ldots &0&0\\
\end{matrix}
\right )
$$
Note that the matrix has two tags: the upper index $(n)$ indicates
the size of the matrix, while the variables enclosed in parentheses
are the total set of variables used in the matrix. This detailed
notation was introduced in \cite{CRS} as several of these matrices
were considered with variable tags throughout. However, here we omit
the list of variables if they are sufficiently clear from the
context.

This matrix will be called a {\sc generic sub--Hankel matrix};
more precisely, $M^{(n)}$ is the {\sc generic sub--Hankel matrix of
order $n$ on the variables} $x_0,\ldots,x_n$.

Throughout we fix a polynomial ring $R=k[\xx]=k[x_0,\ldots,x_n]$ over a field $k$ of characteristic zero. We will
denote by $f^{(n)}(x_0,\ldots,x_n)$ the determinant of
$M^{(n)}(x_0,\ldots,x_n)$ for any $r\geq 1$, and we set $f^{(0)}=1$.

In \cite{CRS} the main objective was to prove that this determinant
is homaloidal and explain the geometric contents of the
corresponding polar map. Here we turn ourselves to the
algebraic-homological behavior of the ideal $J\subset R$ generated
by its partial derivatives.

A common feature between the two approaches is a systematic use of a
recurrence using the subideal  $J_i\subset J$ generated by the first
$i+1$ partial derivatives further divided by the gcd of these
derivatives. We need the following results drawn upon \cite{CRS}.
Throughout we set $f:=f^{(n)}$.

\begin{Lemma}{\rm (\cite[Lemma 4.2]{CRS})}\label{prolegomena}
If $n\geq 2$, then
\begin{itemize}
\item[{\rm (i)}] For $0\leq i\leq n-1$, one has
\begin{equation}\label{g.c.d.s}
\frac {\partial f } {\partial x_0}, \ldots, \frac{\partial f }
{\partial x_i}  \in k\left[x_{n-i},\ldots,x_{n}\right]
\end{equation}
and the g.c.d. of these partial derivatives is $x_n^{n-i-1}$
\item[{\rm (ii)}] For any $i\,$ in the range $1\leq i\leq n-1$, the following holds:
\begin{equation}\label{basic_linear_relation}
x_n\,\frac{\partial f }{\partial
x_i}=-\sum_{k=0}^{i-1}\frac{2i-k}{i}\,\, x_{r-i+k}\, \frac{\partial
f }{\partial x_k}.
\end{equation}
Moreover,
\begin{equation}\label{perfect_linear_relation}
x_n\, \frac{\partial f}{\partial x_n}=(n-1)x_0\,\frac{\partial f
}{\partial x_0}+ (n-2)x_1\, \frac{\partial f}{\partial x_1}+\cdots +
x_{n-2}\,\frac{\partial f} {\partial x_{n-2}}
\end{equation}
\end{itemize}
\end{Lemma}

\begin{Proposition}{\rm (\cite[Proposition 4.3 and its proof]{CRS})}\label{perfect_ideal}
For every $1\leq i\leq n-1$, the ideal $J_i$ is perfect and linearly
presented with recurrent Hilbert--Burch matrix of the form
$$\varphi(J_i)=\left(
\begin{array}{c@{\quad\vrule\quad}c}
\raise5pt\hbox{$2\, x_{n-i}$}&\\
\raise5pt\hbox{${\frac{2i-1}{i}}\, x_{n-i+1}$}&\\
\vdots&\varphi(J_{i-1})\\
\raise5pt\hbox{$\frac{i+1}{i}\, x_{n-1}$}&\\
\multispan2\hrulefill\\
\lower3pt\hbox{$x_n$}&\lower3pt\hbox{{\boldmath $0$}}
\end{array}
\right)
$$
\end{Proposition}

\smallskip

We next prove a few additional results not obtained in \cite{CRS}.

\begin{Lemma}\label{multiplicities_annihilators}
Keeping the above notation, set further $P=(x_{n-1},x_n)\subset
R=k[\xx]$. Then
\begin{enumerate}
\item[{\rm (i)}] $P$ is the radical of $J$ and all the
ideals $J_i, \, 1\leq i\leq n-1$ are $P$-primary.
\item[{\rm (ii)}] The $R_P$-module $R_P/(J_{i})_P$
has length ${i+1} \choose {2}$.
\item[{\rm (iii)}] The ideal of $R_P/(J_{n-1})_P$ generated by
the forms $x_n,x_{n-1}^{n-1}$ has length
${n-1} \choose {2}$.
\end{enumerate}
\end{Lemma}
\demo (i) This is clear from the form of these ideals: any prime
ideal containing any of these has to contain $P$, which is clearly
the unique minimal prime thereof. By
Proposition~\ref{perfect_ideal}, each $J_i$ is perfect, hence
$R/J_i$ is Cohen--Macaulay, thus implying that that $P$ is the only
associated prime of $J_i$.

(ii) By the primary case of the associativity formula for the
multiplicities, one has
$$e(R/J_{i})=\ell(R_P/(J_{i})_P)\,
e(R/P)=\ell(R_P/(J_{i})_P),$$
 since $P$ is generated by linear forms.
On the other hand, from Proposition~\ref{perfect_ideal} one has the
graded free resolution
\begin{equation}\label{perfect_recurrence}
0\rar R(-(i+1))^{i}\rar R(-i)^{i+1}\rar R\rar R/J_{i}\rar 0
\end{equation}
Applying the multiplicity formula for Cohen--Macaulay rings with
pure resolution (\cite{HuMi}), one derives in this case
$e(k[\xx]/J_{i})={{i+1} \choose {2}}$, as required.

(iii) As pointed out previously, one has $(x_n, J_{n-1})=
(x_n,x_{n-1}^{n-1})$. On the other hand,
\begin{eqnarray}
(x_n, J_{n-1})/J_{n-1}&\simeq &(x_n)/(x_n)\cap J_{n-1}\simeq
\left(R/(J_{n-1}:x_n)\right)(1)\\
&=& \left(R/J_{n-2}\right) (1),
\end{eqnarray}
where the equality $(J_{n-1}:x_n)=J_{n-2}$ follows from
Proposition~\ref{perfect_ideal}. Therefore
$$\ell\left((x_n,x_{n-1}^{n-1})_P/(J_{n-1})_P\right)=\ell(R/J_{n-2})_P={{n-1}
\choose {2}},$$
 by part (ii). \qed

\begin{Lemma}\label{equality} With the same notation, one has
$$J/J_{n-1}\simeq \frac{R}{(x_n,x_{n-1}^{n-1})}\,(-(n-1)).$$
\end{Lemma}
\demo
The following isomorphisms of
$R$-graded modules are immediate:
\begin{equation}\label{colonisomorphism}
\frac{J}{J_{n-1}}= \frac{\left(J_{n-1},\frac{\partial f}{\partial
x_n}\right)}{J_{n-1}}\simeq \frac{\left(\frac{\partial f}{\partial
x_n}\right)}{J_{n-1}\bigcap\left(\frac{\partial f}{\partial
x_n}\right) }\simeq \frac{R}{\left(J_{n-1}\colon \frac{\partial
f}{\partial x_n}\right)}\left(n-1\right).
\end{equation}

We claim that $\left(J_{n-1}\colon \frac{\partial f}{\partial
x_n}\right)=(x_n, J_{n-1})$. Once this is proved, we will have $\left(J_{n-1}\colon
\frac{\partial f}{\partial x_n}\right)=(x_n,x_{n-1}^{n-1})$ because
it follows easily from the structure of $f$ and its derivatives that
$(x_n, J_{n-1})=(x_n,x_{n-1}^{n-1})$.

Now, by (\ref{perfect_linear_relation}), $(x_n, J_{n-1})\subset
\left(J_{n-1}\colon \frac{\partial f}{\partial x_n}\right)$ as
trivially $J_{n-1}\subset \left(J_{n-1}\colon \frac{\partial
f}{\partial x_n}\right)$.

For the reverse inclusion,  we proceed as follows.

 Let $r\in (J_{n-1}:\partial f/\partial x_n)$ and write $P:=(x_{n-1},x_n)$.
 Since   $J_{n-1}$ is a $P$-primary ideal and  $\partial f/\partial x_n \not\in J_{n-1} $ then $r\in P=(x_n,x_{n-1})$.
 Next rewrite $r=r(x_0,\ldots,x_n)$ as
 $$r=x_n h(x_0,\ldots,x_n)+r'(x_0,\ldots,x_{n-1})\in P,$$
 where $x_n$ divides no term on the second summand.
Then $r'(x_0,\ldots,x_{n-1})\in (x_{n-1})$, so let $l\in \mathbb{N}$ be such that
$r'(x_0,\ldots,x_{n-1})= x^l_{n-1}r''(x_0,\ldots,x_{n-1})$ and $x_{n-1}$ does not divide $r''(x_0,\ldots,x_{n-1})$.
It follows that $r''(x_0,\ldots,x_{n-1})\not\in P$.

 Recall that $(x_n,J_{n-1})=(x_n,x_{n-1}^{n-1})$. Since $J_{n-1}$ is P--primary and
 $r''x^l_{n-1}\partial f/\partial x_n \in J_{n-1}$ then $x^l_{n-1}\partial f/\partial x_n \in J_{n-1}$.
Thus, for the required reverse inclusion it is enough to show that $l\geq n-1$.
Write
\begin{equation}\label{formula 4}
 x^l_{n-1}\frac{\partial f}{\partial x_n} = \sum^{n-1}_{i=0} r_i\frac{\partial f}{\partial x_i},
\end{equation}
 where $r_i\in R.$

Writing $r_i=x_ng_{i}(x_0,\ldots,x_n)+g'_{i}(x_0,\ldots,x_{n-1})$ and drawing upon (\ref{basic_linear_relation}) yields
 \begin{eqnarray*}
 x^l_{n-1}\partial f/\partial x_n &=& \sum^{n-1}_{i=0} r_i\partial f/\partial x_i
 =\sum^{n-1}_{i=0}(x_ng_{i}(x_0,\ldots,x_n)+g'_{i}(x_0,\ldots,x_{n-1}))\partial f/\partial x_i\\
&=&\sum^{n-1}_{i=0}g_{i}(x_0,\ldots,x_n)(-\sum^{i-1}_{k=0}{ \frac{2i-k}{i}x_{n-i+k} \frac{\partial f}{\partial x_k})
+\sum^{n-1}_{i=0}g'_{i}(x_0,\ldots,x_{n-1})\partial f/\partial x_i}.
\end{eqnarray*}
By repeating this  process for $g_i$ and so forth, after finitly many steps we may suppose that the coefficients of
$\partial f/\partial x_i $ in (\ref{formula 4}) do not involve $x_n$.

Multiplying both sides of (\ref{formula 4}) by $x_n$ and using (\ref{perfect_linear_relation}) yields a syzygy of the
ideal $J_{n-1}$:
\begin{equation}\label{partial}
((n-1)x^l_{n-1}x_0 -x_n r_0)\frac{\partial f} {\partial x_0} + \cdots +(x_{n-2} x^l_{n-1}-x_n r_{n-2})\frac{\partial}{\partial x_ {n-2}}
 - x_n r_{n-1}\frac{\partial f}{\partial x_{n-1}}=0.
\end{equation}

Thinking of this syzygy as a column vector $K$,
we can write
\begin{equation}\label{c}
 K=\alpha_1C_1+\cdots+\alpha_{n-1}C_{n-1}
\end{equation}
for suitable $\alpha_i\in k[x_0,\cdots,x_n]$, where $C_i$ denotes the  $i$th column of the syzygy matrix
of $J_{n-1}$ as in Proposition~\ref{perfect_ideal}:
 $$\begin{pmatrix}
 2x_1 & 2x_2 &2x_3 &\dots &2x_{n-1}\\
\frac {2n-3}{n-1}x_2 &\frac  {2n-5}{n-2}x_3& \frac{2n-7}{n-3}x_4 &\dots& x_n\\
\vdots & \vdots & \vdots & \dots &\vdots\\
\frac {n+2}{n-1}x_{n-3}& \frac{n}{n-2}x_{n-2}&\frac{n-2}{n-3}x_{n-1}&\dots & 0\\
\frac{n+1}{n-1}x_{n-2} & \frac{n-1}{n-2}x_{n-1} &x_n& \dots& 0\\
\frac{n}{n-1}x_{n-1}& x_n &0 & \dots & 0\\
x_n & 0 & 0 & \dots & 0
\end{pmatrix}_{n\times (n-1)}
$$
  This affords the following relations by looking at the last two rows:
$$ x_{n-2} x^l_{n-1}-x_n r_{n-2}=-r_{n-1}\frac{n}{n-1}x_{n-1} + \alpha_2 x_n \quad \mbox{\rm and}
\quad \alpha_1=-r_{n-1}.$$

Since we already assumed $r_i\in k[x_0,\ldots,x_{n-1}]$ for all $i$, then $\alpha_2=-r_{n-2}$.  Therefore
 $$-r_{n-1}\frac{n}{n-1}x_{n-1}=x_{n-2} x^l_{n-1}.$$
 In particular $l\geq 1$ and hence $ r_{n-1}=-\frac{n-1}{n} x_{n-2} x^{l-1}_{n-1}.$

 \smallskip

 Inspecting the $(n-2)$'th row yields:
 $$2x^l_{n-1}x_{n-3}-x_n r_{n-3}=-\frac{n+1}{n-1} r_{n-1}x_{n-2}-\frac{n-1}{n-2}r_{n-2}x_{n-1}+\alpha_3 x_n.$$
 Since $r_{n-1},r_{n-2}\in k[x_0,\ldots,x_{n-1}]$ then $\alpha_3=-r_{n-3}$.
 Substituting for $r_{n-1}$ obtains
 \begin{eqnarray*}
 2x^l_{n-1}x_{n-3}&=& -\frac{n+1}{n-1} r_{n-1}x_{n-2}-\frac{n-1}{n-2}r_{n-2}x_{n-1}\\
&=& \frac{n+1}{n} x^2_{n-2} x^{l-1}_{n-1}-\frac{n-1}{n-2}r_{n-2}x_{n-1}
\end{eqnarray*}
Then necessarily $l\geq 2$  and furthermore
 \begin{equation}\label{sn}
   r_{n-2}=\frac {n-2}{n-1}(-2x^{l-1}_{n-1}x_{n-3}+\frac{n+1}{n}  x^2_{n-2} x^{l-2}_{n-1})=x^{l-2}_{n-1}s_{n-2},
 \end{equation}
 $\text{for some }  s_{n-i}\in k[x_0,\ldots,x_{n-1}].$

  Since $x_n$ appears exactly once on each row of the syzygy matrix below the first one, the argument inducts yielding
  that for every $1\leq i \leq n-2$ one has $\alpha_i=-r_{n-i}$ and $r_{n-i} = x^{l-i}_{n-1}s_{n-i}$
  with $s_{n-i}\in k[x_0,\ldots,x_{n-1}]$ and $l\geq i$.
    In particular,  $l\geq n-2$ and  $r_2= x^{l-n+2}_{n-1}s_2$.

\smallskip

 Finally, from the first row of $K$, we have:
 $$(n-1)x_0x^l_{n-1}-x_nr_0=-2x_1r_{n-1}-2x_2r_{n-2}-\ldots-2x_{n-2}r_2-2x_{n-1}r_1.$$
Hence $r_0=0$. Rearranging yields
\begin{equation}\label{r1}
 -2x_1r_{n-1}-2x_2r_{n-2}-\ldots-2x_{n-1}r_1-(n-1)x_0x^l_{n-1}=2x_{n-2}r_2=2x_{n-2}x^{l-n+2}_{n-1}s_2.
 \end{equation}
Since the left hand side is divisible by $x_{n-1}$ so is the right hand side.
 Thus $l-n+2>0$. In other words, $l\geq n-1,$ as desired.
\qed

\begin{Proposition}\label{resolution of J}
Let $J$ denote the gradient ideal of the sub-Hankel determinant. Then the minimal graded resolution of $R/J$ has the form
\begin{equation*}
0\rar R(-(2n-1))\rar R(-n)^n\oplus
R(-2(n-1))\stackrel{\phi}{\rar} R^{n+1}(-(n-1))\rar R.
\end{equation*}
In particular, $R/J$ is strictly almost Cohen--Macaulay.
\end{Proposition}
\demo
By Lemma~\ref{equality} we have a free minimal resolution
\begin{equation*}
{\cal C}\colon\quad 0\rar R(-(2n-1))\rar R(-n)\oplus R(-2(n-1))\rar R(-(n-1))
\rar
J/J_{n-1}\rar 0,
\end{equation*}
On the other hand, (\ref{perfect_recurrence}) gives a resolution
\begin{equation*}
{\cal J}_{n-1}\colon\quad 0\rar R(-n)^{n-1}\rar R(-(n-1))^n\rar
R\rar R/J_{n-1}\rar 0
\end{equation*}

 Since the inclusion $J/J_{n-1}\subset R/J_{n-1}$ induces a map of complexes ${\cal C}\rar {\cal J}_{n-1}$,
the resulting mapping cone is a resolution of $R/J$ (\cite[Exercise A3.30]{Eisenbook}):
\begin{equation}\label{resolution3}
\quad 0\rar R(-(2n-1))\rar R(-n)^n\oplus
R(-2(n-1))\stackrel{\phi}{\rar} R^{n+1}(-(n-1))\rar R\rar R/J\rar 0,
\end{equation}
(where the right end tail $R\oplus J/J_{n-1}\rar R/J_{n-1}\rar 0$
has been replaced by $R\rar R/J\rar 0$).
Moreover, since for every relevant index $i$, the shifts of $({\cal
J}_{n-1})_i$ are strictly smaller than those of $({\cal C})_i$, it
follows by [loc. cit.] that (\ref{resolution3}) is minimal. In
particular, one reads from it that the linear part $\phi_1$ has rank
$n$, hence maximal.
\qed

\medskip

Let us dwell a little more on the details of the mapping cone in the previous proof.
It is of the form

$$\begin{CD}
0 @>>> R^{n-1} @> \varphi_{n-1} >> R^n @> (\frac{\partial f}{\partial x_0},\ldots,\frac{\partial f}{\partial x_{n-1}})
>> R@>>> \frac {R}{J_{n-1}}@>>> 0\\
&&   @A {g_4}AA   @A{g_3}AA    @A{g_2}AA @A{g_1} AA && \\
 0 @>>>       R @> (x^{n-1}_{n-1},-x_n)^t>>    R^2  @>(x_n,x^{n-1}_{n-1})>>   R @>>>     J/J{n-1} @>>> 0
\end{CD}\
$$

\smallskip

Note that $g_2$ is  multiplication by $\partial f/ \partial x_n$ and the induced map $g_3$ is given by the following $n\times2$ matrix:

 $$g_3=\begin{pmatrix}
 (n-1)x_0 & r_0\\
 (n-2)x_1 & r_1\\
 \vdots\\
 x_{n-2} & r_{n-2}\\
 0 & r_{n-1}
 \end{pmatrix}.$$

We next find out the entries of $g_4$. By the commutativity of  the mapping cone diagram, we have
$$\varphi_{n-1}\circ g_4=g_3 \circ \begin{pmatrix}
x^{n-1}_{n-1}\\
-x_n
\end{pmatrix}= \begin{pmatrix}
 (n-1)x_0 & r_0\\
 (n-2)x_1 & r_1\\
 \vdots\\
 x_{n-2} & r_{n-2}\\
 0 & r_{n-1}
 \end{pmatrix}\begin{pmatrix}
x^{n-1}_{n-1}\\
-x_n
\end{pmatrix}=\begin{pmatrix}
(n-1)x^{n-1}_{n-1}x_0 -x_n r_0 \\
(n-2)x^{n-1}_{n-1}x_1 - x_n r_1 \\
\vdots\\
x_{n-2} x^{n-1}_{n-1}-x_n r_{n-2}\\
-x_n r_{n-1}
\end{pmatrix} $$
where the rightmost matrix is the syzygy $K$ in (\ref{partial}), viewed as a column vector, where $l=n-1$.
By reasoning as in the argument that ensues (\ref{c}), one gets that the $i$th entry of $g_4$ is $- r_{n-i}$, $\,1\leq i\leq n-1$.
As a result, the leftmost map in (\ref{resolution3}) is  $\psi:=(x^{n-1}_{n-1}, -x_n, g_4^t)=(x^{n-1}_{n-1}, -x_n,-r_{n-1},\ldots,-r_1)$.

\smallskip

We make use of this in the following result, where $J$ is the gradient ideal of the sub-Hankel determinant.

\begin{Theorem}
The associated primes of $R/J$ are $(x_{n-1}, x_n)$ and $(x_{n-2},x_{n-1}, x_n)$.
\end{Theorem}
\demo
Clearly, $P:=(x_{n-1}, x_n)$ is a minimal prime thereof.
Since $P$ is also the radical of $J$ (Proposition~\ref{multiplicities_annihilators} (i))
then there are no additional minimal primes.

To argue for embedded primes we proceed as follows.
Let as above $\psi$ denote the tail map of (\ref{resolution3}).
Since $R/J$ has homological dimension $3$,  any $Q\in \Ass (R/J)$ has codimension at most $3$ and a prime $Q$ of
codimension $3$ containing $J$ is an associated prime of $R/J$
if and only if $Q\supset I_1(\psi)$ (see, e.g., \cite[Corollary 20.14(a)]{Eisenbook} for the last part).
As seen above, $I_1(\psi)=(x^{n-1}_{n-1}, x_n,-r_{n-1},\ldots,-r_1)$.
On the other hand, we know
from the proof of Proposition~\ref{equality} that  $r_{n-i} = x_{n-1}^{n-i-1}s_{n-i}$,
where $s_{n-i}\in k[x_0,\ldots,x_{n-1}]$ for all $1\leq i \leq n-1$.
Additionally, the inductive argument in this proof and relation (\ref{sn}) also show that, for $1\leq i \leq n-2$,
$s_{n-i}=x_{n-1}t_{n-i}+x_{n-i}^i $,  where $t_{n-i}\in  k[x_0,\ldots,x_{n-1}]$ .

 Therefore the condition for such a prime $Q$ to be an associated prime of $R/J$ is that it contain the ideal
 $(x_n, x_{n-1},x_{n-2}s_2)=(x_n, x_{n-1},x_{n-2}(x_{n-1}t_2+x_{n-2}^{n-2}))=(x_n, x_{n-1},x_{n-2}^{n-1} )$.
 It follows that $Q=(x_{n-2},x_{n-1}, x_n)$, as stated.
\qed

\begin{Corollary}
Set $P:=\sqrt{J}=(x_{n-1}, x_n)$.
Then the $P$-primary component of $J$ is $J_{n-2}$.
\end{Corollary}
\demo
Since $J_{n-2}$ is a $P$-primary ideal (Proposition~\ref{multiplicities_annihilators} (i)), it is equivalent to show the equality
$J_P={(J_{n-2})}_P$.  Since $J\subset J_{n-2}$ we will be done by showing the equality of lengths
$\emph{l}(\frac {R_P}{J_P})=\emph{l}( \frac {R_P}{{(J_{n-2})}_P})$.

Now, with the present data, by the associativity formula one has $\emph{l}(\frac {R_P}{J_P})=e(R/J)$ -- the multiplicity
of $R/J$.
To compute the latter we deploy the numerator of the Hilbert series of $R/J$ in terms of the graded Betti numbers of $R/J$
as in  (\ref{resolution3}); it obtains
$$S(t):= 1-(n+1)t^{n-1} + t^{2n-2} + nt^n - t^{2n-1}.$$
Taking second derivatives, evaluating at $t=1$, etc., finally gives
 $e(R/J)= (n-1)(n-2)/2 = {n-1\choose 2}$. But the latter coincides with  $\emph{l}(R_P / {(J_{n-2})}_P)$ by
 Lemma~\ref{multiplicities_annihilators} (ii).
\qed

 \begin{Theorem} The gradient ideal $J$ of the sub-Hankel determinant is of linear type.
 \end{Theorem}
 \demo
 By definition, we have to show that the natural surjective $R$-homomorphism $\mathcal{S}_R(J)\surjects \mathcal{R}_R(J)$ from
the symmetric algebra of $J$ to its Rees algebra is injective.
One knows that this is the case if and only if $\mathcal{S}_R(J)$ is a domain.

Set $\mathcal{S}_R(J)\simeq R[y_0,\cdots,y_n]/\mathcal{L}$, where $\mathcal{L}$ is the ideal generated by the $1$-forms coming from
the syzygies of $J$.
We will argue as follows: since $x_n$  belongs to the radical of $J$,  $J_{x_n}$ is the unit ideal in $R_{x_n}$.
Now, suppose one shows that $x_n$ is a non zero-divisor modulo $\mathcal{L}$. Then  $\mathcal{L}_{x_n}$ is
the defining ideal of the symmetric algebra of $J_{x_n}=R_{x_n}$, hence it is the zero ideal in a polynomial ring over a domain.
In particular, it is a prime ideal, hence so must be $\mathcal{L}$. Therefore $\mathcal{S}_R(J)$ is a domain, thus
the structural homomorphism is injective as observed.

By a quirk, it will be easier to show first that $y_n$ is non-zero divisor modulo $\mathcal{L}$ and then that
 $(y_n,x_n)$ is a regular sequence modulo $\mathcal{L}$. In this case, since $\mathcal{S}_R(J)$ is a
 positively graded ring any permutation of a regular sequence is a regular sequence, hence $(x_n,y_n)$ is a regular sequence
 as well, in particular $x_n$ is a non-zero divisor over $\mathcal{S}_R(J)$.

\textbf{Step 1.} $y_n$ is non-zero divisor modulo $\mathcal{L}$.

 Let $h=h(\textbf{x},\textbf{y})\in R[y_0,\ldots,y_n]$ be such that $y_n h(\textbf{x},\textbf{y})\in \mathcal{L} $.
 Say,
 $$y_n h=\sum_{i=1}^{n+1}h_i g_i,$$
 for suitable $h_i=h_i(\xx,\yy) \in R[y_0,\ldots,y_n]$, where
$$g_i=2x_{n-i}y_0 +\frac{2i-1}{i}x_{n-i+1}y_1+ \cdots +x_ny_i  \;\;(1\leq i\leq n-1)$$
 and
 $$g_n = (n-1)x_0y_0+(n-2)x_1y_1+\cdots +x_ny_n,\;\; g_{n+1}=\sum_{j=1}^{n-1} {r_jy_j} +x_{n-1}^{n-1}y_n$$
 generate $\mathcal{L}$ (from (\ref{resolution3})), with $r_i$ being as in (\ref{partial}).

 Decompose further $h_i={h'_i}+y_n{h''_i}$, with ${h'_i}\in  R[y_0,\ldots,y_{n-1}]. $
 Then
 $$y_n h(\textbf{x},\textbf{y})=y_n\sum_{i=1}^{n-1}{h''_i}g_i + y_n{h''_n}g_n + {h'_n}x_ny_n + y_n{h''_{n+1}}g_{n+1} + {h'_{n+1}}x_{n-1}^{n-1}y_n.$$
 Since in the right hand side the terms not divisible by $y_n$ must vanish, we get
  $$h=\biggl (\sum_{i=1}^{n-1}{h''_i} g_i + {h''_n}g_n + {h''_{n+1}}g_{n+1}\biggr) + {h'_n}x_n + {h'_{n+1}}x_{n-1}^{n-1}.$$
 To show that $h\in \mathcal{L}$ it is then enough to check that
 $$ {h'_n}x_n+ {h'_{n+1}}x_{n-1}^{n-1}\in\mathcal{L}. $$
 Now a form in $\mathcal{L}$ must vanish when evaluated at $y_i\mapsto \partial f/\partial x_i$ -- the generators of $J$.
 Letting $\partial\ff$ denote these partial derivatives we have
  $\frac{\partial f}{\partial x_n}h(\textbf{x},\partial\ff)=0$, hence $h(\textbf{x},\partial\ff)=0.$
  Retrieving in terms of the expression of $h$ implies that
 $$ h'_n(\xx,\partial\ff)\,x_n+ h'_{n+1}(\xx,\partial\ff)\,x_{n-1}^{n-1}=0.$$
 Since $h'_n ,h'_{n+1}\in R[y_0,\ldots,y_{n-1}]$, the form $h'_nx_n+ h'_{n+1}x_{n-1}^{n-1}$ belongs to the defining
 ideal of $\mathcal{R}_R(J_{n-1})$.
 By \cite[Proposition 4.3]{CRS}, $J_{n-1}$ is of linear type, hence  $h'_nx_n+ h'_{n+1}x_{n-1}^{n-1}$ belongs to
 the defining ideal of $\mathcal{S}_R(J_{n-1})$, which is a subideal of $\mathcal{L}$ by the general theory
 developed in \cite{CRS}.
 This shows the contention.

 \medskip

 \textbf{Step 2.} $x_n$ is non-zero divisor modulo $(\mathcal{L},y_n)$.

 \smallskip

One has $(\mathcal{L},y_n)=(g_1,\ldots,g_{n-1},g,h,y_n),$ where
$g=(n-1)x_0y_0+(n-2)x_1y_1+\cdots+x_{n-2}y_{n-2}$ and $h=\sum_{j=1}^{n-1} {r_jy_j}.$
 Then
  $$\frac{R[y_0,\ldots,y_n]}{(\mathcal{L},y_n)}\simeq
  \frac{R[y_0,\ldots,y_{n-1}]}{(g_1,\ldots,g_{n-1},g,h)}$$
hence we are to show that $x_n$ is a nonzerodivisor on the rightmost ring.
  Let then $\kappa=\kappa(\mathbf{x}, \mathbf{y}\setminus y_n)$ be a form in $R[y_0,\ldots,y_{n-1}]$
  such that $x_n\kappa\in (g_1,\ldots,g_{n-1},g,h)$.
  Write $x_n\kappa=\sum_{j=1}^{n-1} \mu_jg_j+\mu g +\nu h$, for suitable forms $\mu_j,\mu,\nu\in R[y_0,\ldots,y_{n-1}]$.
 Evaluating $y_i\mapsto f_{x_i}:=\partial f/\partial x_i$ for  $i=0,\cdots,n-1$, and taking in account
 the shape of $g$ and $h$, we have
 \begin{eqnarray*}
 x_n\kappa(\mathbf{x},f_{x_0},\ldots, f_{x_{n-1}})&=&\mu(\mathbf{x},f_{x_0},\ldots, f_{x_{n-1}})g(\mathbf{x},f_{x_0},\ldots, f_{x_{n-1}})\\
      &+&\nu(\mathbf{x},f_{x_0},\ldots, f_{x_{n-1}})h(\mathbf{x},f_{x_0},\ldots, f_{x_{n-1}})\\
      &=& -\mu(\mathbf{x},f_{x_0},\ldots, f_{x_{n-1}})x_n\,f_{x_n} -\nu(\mathbf{x},f_{x_0},\ldots, f_{x_{n-1}})x_{n-1}^{n-1}\,f_{x_n}.
      \end{eqnarray*}
On the other hand, by the shape of $f_{x_n}$,  one has $\gcd(x_n, f_{x_n})=1$.
It follows that
$$\kappa(\mathbf{x},f_{x_0},\ldots, f_{x_{n-1}})=\delta\,f_{x_n},$$
for some $\delta\in R$.
Pulling back to the $\mathbf{y}$-variables tells us that $\kappa  -\delta\, y_n$
vanishes on the partial derivatives, and hence it belongs to affine ideal $\tilde{\mathcal{J}}$
of all polynomials vanishing on the partial derivatives.
This ideal is prime because we can consider $\partial\ff:=(f_{x_0},\ldots , f_{x_n})$ as a point in $K^{n+1}$, where $K$ denotes the field
of fractions of $R$ and consider the ideal of $K[\mathbf{y}]$ vanishing on $\partial\ff$ and then contract to $R[\mathbf{y}]$.
We note that the Rees ideal $\mathcal{J}$ is the largest homogeneous ideal contained in $\tilde{\mathcal{J}}$.

Now, multiplying $\kappa  -\delta\, y_n$ by $x_n$ and using that $x_n\kappa\in \mathcal{L}\subset \mathcal{J}
\subset \tilde{\mathcal{J}}$, it follows that $x_n \,\delta\, y_n\in  \tilde{\mathcal{J}}$.
Since this element is (trivially) homogeneous in $\mathbf{y}$, it must belong to the Rees ideal, hence $\delta=0$.
Therefore, $\kappa\in \tilde{\mathcal{J}}$. But, since $\kappa$ is assumed to be homogeneous, it belongs to the Rees ideal.
Thus, we have $\kappa\in \mathcal{J}\cap R[y_0,\ldots, y_{n-1}]$.
But this means that $\kappa$ belongs to the Rees ideal $\mathcal{J}'$ of $J_{n-1}$.
Since the latter is of linear type by we conclude as above that $\kappa\in \mathcal{L}$, as was to be shown.
\qed

\subsection{Degenerations of a catalecticant}

The moral of this section is to search for predecessors of the sub-Hankel matrix, coming down all the way from the generic matrix.
The case of the sub-Hankel determinant taught us that degenerating (with zeros) improves homaloidness and actually
also the good algebraic properties of the polar map (such as the number of linear syzygies thereof).
Guided by this case, we may expect a similar situation while taking degenerations of the generic, more generally, catalecticants.
In what follows we will see that our hopes of analogy are a bit naive.

\subsubsection{Degenerating the generic matrix}\label{deggeneric}

We consider the simplest degeneration of the $m\times m$ generic matrix, consisting of replacing a variable by zero.
Since in the generic matrix an entry has equal share as any other entry, we may assume that the resulting matrix has the shape:

\begin{equation}\label{degener_gen}
\mathcal{DG}=\mathcal{DG}_m:=\left(
  \begin{array}{cccc|c}
    x_0 & x_1 & \ldots & x_{m-2} & \mathbf {x_{m-1}} \\
    x_m & x_{m+1} & \ldots & x_{2m-2}& \mathbf{x_{2m-1}} \\
    \vdots & \vdots & \vdots & \vdots & \vdots\\
    x_{m(m-2)} & x_{m(m-2)+1} & \ldots & x_{m(m-1)-2} & \mathbf{x_{m(m-1)-1}}\\[5pt]
    \hline\\[-12pt]
    \mathbf{x_{m(m-1)}} &  \mathbf{x_{m(m-1)+1}}& \ldots & \mathbf{x_{m^2-2}} & 0 \\
  \end{array}
\right)
\end{equation}
The reason for the boldfaced variables in a minute.
We assume that $m\geq 3$.

\begin{Proposition}
Let $R=k[x_0,\ldots , x_{m^2-2}]$, let $f:=\det\mathcal{DG}\in R$ and let $J\subset R$ denote the gradient ideal of $f$.
Then:
\begin{enumerate}
\item[{\rm (a)}] $f$ is a  Gordan--Noether polynomial{\rm ;} in particular, its Hessian determinant vanishes.
\item[{\rm (b)}] The image of the polar map of $f$ is the locus of the  submaximal minors of the $(m-1)\times (m-1)$ square submatrix to the left and above the crossing lines{\rm ;} in other words, it is the cone over the Segre embedding $\pp^{m-2}\times \pp^{m-2}$.
\item[{\rm (c)}] $J$ has maximal linear rank and the associated primes of $R/J$ are the Gorenstein ideals $I_{m-1}(\mathcal{DG})$ and the ideal generated by the boldfaced entries in {\rm (\ref{degener_gen})}.
\end{enumerate}
\end{Proposition}
\demo
(a) Expanding $f$ by Laplace along the last row, then expanding the corresponding (signed) cofactors along the last columns yields the following expression of $f$:
\begin{equation}
f=\bold\Delta\cdot \bold\Sigma ^T,
\end{equation}
where $\bold\Delta$ denotes the row of the signed $(m-2)$-minors of the $(m-1)\times (m-1)$ generic square submatrix $\mathcal{G}_{m-1}$ to the left and above the crossing lines in (\ref{degener_gen}), while $\bold\Sigma ^T$ stands for the transpose of the row vector whose entries are all the cross products of the entries of the last column and the last row (not counting the zero at the corner).

Now, since the entries of $\bold\Sigma ^T$ are algebraically dependent if $m\geq 3$ and they share no variables with the entries of  $\bold\Delta$ then $f$ is a particular case of a Gordan--Noether polynomial (see \cite[Section 2.3]{CRS}).
Therefore, $f$ has vanishing Hessian.

(b) Since $f$ has the structure of a Gordan--Noether polynomial, one knows that the defining equations of the polar variety (i.e., the image of the polar map) of $f$ are exactly the polynomial relations of the entries in  $\bold\Sigma ^T$. But the latter are the parameters defining the Segre embedding of $\pp^{m-2}\times \pp^{m-2}$, whose
defining equations are exactly the $(m-2)$-minors of the generic matrix  $\mathcal{G}_{m-1}$ read in the dual variables.

(c) We observe that the ideal  $I_{m-1}(\mathcal{DG})$ is a  Gorenstein ideal of codimension $4$ since it specializes from the full generic case, and as such it is a linearly presented ideal.
Moreover, it is clear that since the gradient ideal of the fully generic predecessor of $\mathcal{DG}$ is its ideal of $(m-1)$-minors then one has $I_{m-1}(\mathcal{DG})=(J, \det \mathcal{G}_{m-1})$.
From this follows easily that the number of $k$-linearly independent syzygies of $J$ falls by just a bit off the ones of  $I_{m-1}(\mathcal{DG})$. But the latter is linearly presented with a number of syzygies $>> m^2-2$. Thus, $J$ has maximal linear rank $m^2-2$.

By a codimension argument, $I_{m-1}(\mathcal{DG})$ is a minimal prime of $R/J$ (actually, it is the unmixed part of the gradient ideal).
By the shape stressed in (\ref{degener_gen}) it is clear that the boldfaced variables conduct $\det \mathcal{G}_{m-1}$ into $J$ (for this use the Cramer--Hilbert--Burch relations along the obvious $(m-1)\times m$ and $m\times (m-1)$ submatrices of  $\mathcal{DG}$).
The other inclusion is easy as well.
Therefore,
$$J:I_{m-1}(\mathcal{DG})=J:\det \mathcal{G}_{m-1}=(\mbox{\rm boldfaced entries}),$$
which implies the assertion.
\qed

\subsubsection{Degeneration of a $2$-leap catalecticant}

Consider the following $3\times 3$ matrix, obtained from the $2$-leap generic catalecticant by
replacing the last entry with zero:

$$\mathcal{SC}_3=\left(
  \begin{array}{ccc}
    x_0 & x_1 & x_2 \\
    x_2 & x_3 & x_4 \\
    x_4 & x_5 & 0 \\
  \end{array}
\right)
$$
\begin{Proposition}
Let $f:=\det \mathcal{SC}_3$ and let $J$ denote the gradient ideal of $f$.
Then:
\begin{enumerate}
\item[{\rm (a)}] $f$ has non-vanishing Hessian determinant
\item[{\rm (b)}] $J$ has maximal linear rank.
\end{enumerate}
In particular, $f$ is homaloidal.
\end{Proposition}
\demo For the proof we really need to resort to the computer:
we cannot use Proposition~\ref{2-cat3x3} as guidance since even the ideal of $2$-minors does not specialize, its codimension now being $3$ and not Gorenstein (not even Cohen--Macaulay).

Still, the shape of the Hessian matrix is nearly subdiagonal with $x_4$ along the main subdiagonal -- a situation pretty much like the one in the sub-Hankel case. This makes it easy to inspect the value of the Hessian determinant to conclude that it is indeed a power of $x_4$ up to a nonzero coefficient.
This takes care of (a).

For (b), we do not learn immediately from the syzygies of the $2$-minors even knowing that the latter are again linearly presented -- the reason being that only four of the partial derivatives are $2$-minors.
Resorting to \cite{Macaulay} gives that $J$ has $7$ linear syzygies generating a submodule of rank $5=6-1$.
\qed

\begin{Remark}\rm
Collecting the results so far, one has seen three sorts of behavior of a determinant subjected to degeneration: (1) It may not be homaloidal before degenerating and become homaloidal after (sub-Hankel); (2) it may be homaloidal before degeneration and stop being so afterwards (sub-Generic); and (3) it can be homaloidal before and continue being homaloidal after degeneration ($2$-leap sub-catalecticant.)
\end{Remark}

\bigskip

\noindent {\bf Authors' address:}

\medskip

\noindent Departamento de Matem\'atica, CCEN\\ Universidade Federal
de Pernambuco\\
50740-560 Recife, PE, Brazil.\\
maral@dmat.ufpe.br,  aron@dmat.ufpe.br

\end{document}